\documentclass[12pt,leqno]{article}
\usepackage[utf8]{inputenc}
\usepackage{tikz}
\usepackage[official]{eurosym}
\usepackage{mathtools,fullpage,amsmath,amssymb,theorem}
\usepackage{accents}

\usepackage{epic,eepic}
\usepackage{cases}

\def\Mfor{\mathbb{M}[\vec{U}]}

\def\MforI{\mathbb{M}_I[\vec{U}]}

\def\MforIW{\mathbb{M}_I[\vec{W}]}
\def\min{{\rm min}}
\def\max{{\rm max}}
\def\sup{{\rm sup}}
\def\otp{{\rm otp}}
\def\Lim{{\rm Lim}}

\def\llvdash{{\|\hskip-2pt \raise 3pt\hbox{\vrule
height 0.25pt width 0.4cm}}}

\def\l{{\langle}}
\def\r{{\rangle}}

\def\oa{{\overline A^{\,\lower 7pt_{\hbox{$\scriptstyle\bet}}
\hbox{$\scriptstyle 0\tau$}}}}

\def\bet{\beta}

\def\llvdash{{\|\hskip-2pt \raise 3pt\hbox{\vrule height 0.25pt
width 0.4cm}}}

\newtheorem{theorem}{Theorem}[section]
\newtheorem{lemma}[theorem]{Lemma}
\newtheorem{corollary}[theorem]{Corollary}
\newtheorem{proposition}[theorem]{Proposition}
{\theorembodyfont{\rmfamily}
\newtheorem{definition}[theorem]{Definition}}
{\theorembodyfont{\rmfamily}
\newtheorem{remark}[theorem]{Remark}}
{\theorembodyfont{\rmfamily}
\newtheorem{example}[theorem]{Example}}
{\theorembodyfont{\rmfamily}
\newtheorem{claim}{Claim}}
{\theorembodyfont{\rmfamily}
}
{\theorembodyfont{\rmfamily}
}

\newcommand{\pr}{\medskip\noindent\textit{Proof}. }

\newcommand{\lusim}[1]{\smash{\underset{\raisebox{1.2pt}[0cm][0cm]{$\sim$}}
{{#1}}}}

\def\otp{{\rm otp}}

\def\llvdash{{\|\hskip-2pt \raise 3pt\hbox{\vrule
height 0.25pt width 0.15cm}}}

\def\Vdashbks{\hbox{$\Vdash\!\!\!\!{\raise2pt\hbox
{$\scriptscriptstyle\backslash$}}$}}

\title{Intermediate Models of Magidor-Radin Forcing-Part I}
\author{ Tom Benhamou and Moti Gitik\footnote{ The work of the second author was partially supported by ISF grant No.1216/18.}}
\date{\today}

\begin{document}

\maketitle
\begin{abstract}
    We continue the work done in \cite{PrikryCaseGitikKanKoe},\cite{TomMoti}. We prove that for every set $A$ in a Magidor-Radin generic extension using a coherent sequence such that $o^{\vec{U}}(\kappa)<\kappa$, there is a subset $C'$ of the Magidor club such that $V[A]=V[C']$. Also we classify all intermediate $ZFC$ transitive models $V\subseteq M\subseteq V[G]$.
\end{abstract}
\section{Introduction}
 In this paper we consider the version of Magidor-Radin forcing for $o^{\vec{U}}(\kappa)\leq\kappa$, but prove results for $o^{\vec{U}}(\kappa)<\kappa$. Section $(2)$, will also be relevant to the forcing in Part II. 
 
 Denote by $C_G$, the generic Magidor-Radin club derived from a generic filter $G$.
In \cite{TomMoti}, the authors proved the following: 
 \begin{theorem}\label{result of master}
 Let $\vec{U}$ be a coherent sequence and $G\subseteq \Mfor$ be a $V$-generic filter such that $o^{\vec{U}}(\beta)<\delta_0:=\min\{\alpha\mid 0<o^{\vec{U}}(\alpha)\}$ for every $\beta\in C_G\cup\{\kappa\}$. Then for every set $A\in V[G]$, there is $C\subseteq C_G$ such that $V[A]=V[C]$.
 \end{theorem} 
 
  In this paper we would like to generalize this result to the case where $o^{\vec{U}}(\kappa)<\kappa$. 
 Formally, we prove this generalization by induction $\kappa$, namely, the inductive hypothesis is that for every $\delta<\kappa$, any coherent sequence $\vec{W}$ with maximal measurable $\delta$, and any set $A$ in a generic extension $V[H]$, where $H\subseteq\mathbb{M}[\vec{W}]$, there is $C\subseteq C_H$ such that $V[A]=V[C]$. Here we do not restrict the order of $\delta$'s below $\kappa$. To be precise, the proof given in this paper is the inductive step for the case $o^{\vec{U}}(\kappa)<\kappa$:
 \begin{theorem}
 Let $U$ be a coherent sequence with maximal measurable $\kappa$, such that $o^{\vec{U}}(\kappa)<\kappa$. Assume the inductive hypothesis that for every $\delta<\kappa$, any coherent sequence $\vec{W}$ with maximal measurable $\delta$, and any set $A$ in a generic extension $V[H]$ for $H\subseteq\mathbb{M}[\vec{W}]$, there is $C\subseteq C_H$ such that $V[A]=V[C]$. Then for every $V$-generic filter $G\subseteq\Mfor$ and any set $A\in V[G]$, there is $C\subseteq C_G$ such that $V[A]=V[C]$.
 \end{theorem}
 As a corollary of this, we obtain the main result of this paper:
 \begin{theorem}\label{MainResaultPartone}
 Let $\vec{U}$ be a coherent sequence such that $o^{\vec{U}}(\kappa)<\kappa$. Then for every $V$-generic filter $G\subseteq\Mfor$, such that $\forall\alpha\in C_G. o^{\vec{U}}(\alpha)<\alpha$ and every $A\in V[G]$, there is $C\subseteq C_G$ such that $V[A]=V[C]$.
 \end{theorem}

The first problem which rises when we let $o^{\vec{U}}(\kappa)\geq\delta_0$, is that we might lose completness for some of the pairs in a condition $p$. For example, if $p=\langle\langle \delta_0,A_0\rangle,\langle\kappa,A_1\rangle\rangle$, we wont be able to take in account all the measures on $\kappa$, since there are $\delta_0$ many of them and only $\delta_0$-completness.  
 The idea is to split $\Mfor$ to the part below $o^{\vec{U}}(\kappa)$ and above it. The cardinality of the lower part is lower than the the degree of $\leq^*$-closure of the upper part. The upper part is an instance of $\Mfor$, where the order of every measurable is below the order of $\kappa$ which is similar to the framework of theorem \ref{result of master}, then some but not all of the arguments of \cite{TomMoti} generalize.

Note that the classification we had in \cite{TomMoti} for models of the form $V[C']$, does not extend, even if $o^{\vec{U}}(\kappa)=\delta_0$. 
\begin{example}\label{examplenonGeneralize}
Consider $C_G$ such that $C_G(\omega)=\delta_0$ and $o^{\vec{U}}(\kappa)=\delta_0$.
Then in $V[G]$ we have the following sequence $C'=\langle C_G(C_G(n))\mid n<\omega\rangle$ of points of the generic $C_G$ which is determine by the first Prikry sequence at $\delta_0$.

Then $I(C',C_G)=\langle C_G(n)\mid n<\omega\rangle\notin V$, where $I(X,Y)$ is the indices of $X\subseteq Y$ in the increasing enumeration of $Y$.

The forcing $\mathbb{M}_I[\vec{U}]$ which was defined in \cite{TomMoti}, is no longer defined in $V$ since $I\notin V$. 

In this case, we will add points to $C'$, which are simply $\langle C_G(n)\mid n<\omega\rangle$, then the forcing will be a two step iteration.
The first will be to add the Prikry sequence $\langle C_G(n)\mid n<\omega\rangle$, then the second will be a Diagonal Prikry forcing adding point from the measures $\langle U(\kappa,C_G(n))\mid n<\omega\rangle$, which is of the form $M_{I}[\vec{U}]$.
\end{example}
Generally, we will define forcing $\mathbb{M}_f[\vec{U}]$, which are not subforcing of $\Mfor$, but are a natural diagonal generalization of $\Mfor$ and a bit closer to Magidor's original formulation in \cite{ChangeCofinality}.

The classification of models is given by the following theorem:
\begin{theorem}\label{Classification}
Assume that for every $\alpha\leq\kappa$, $o^{\vec{U}}(\alpha)<\alpha$. Then for every $V$-generic filter $G\subseteq\Mfor$ and every transitive $ZFC$ intermediate model $V\subseteq M\subseteq V[G]$, there is a closed subset $C_{fin}\subseteq C_G$ such that:
\begin{enumerate}
    \item $M=V[C_{fin}]$.
    \item  There is a finite iteration $\mathbb{M}_{f_1}[\vec{U}]*\mathbb{M}_{\underaccent{\sim}{f}_2}[\vec{U}]...*\mathbb{M}_{\underaccent{\sim}{f}_{n}}[\vec{U}]$, and a $V$-generic $H^*$ filter for $\mathbb{M}_{f_1}[\vec{U}]*\mathbb{M}_{\underaccent{\sim}{f}_2}[\vec{U}]...*\mathbb{M}_{\underaccent{\sim}{f}_n}[\vec{U}]$ such that $V[H^*]=V[C_{fin}]=M$.
\end{enumerate}
\end{theorem}
\section{ Basic Definitions and Preliminaries}
We will follow the description of Magidor forcing as presented in \cite{Gitik2010}. 
 
 Let $\vec{U}=\langle U(\alpha,\beta)\mid \alpha\leq \kappa \ ,\beta<o^{\vec{U}}(\alpha)\rangle$ be a coherent sequence. For every $\alpha\leq\kappa$, denote $$\cap\vec{U}(\alpha)=\underset{i<o^{\vec{U}}(\alpha)}{\bigcap}U(\alpha,i)$$
\begin{definition}\label{Magidor-conditions}
$\mathbb{M}[\vec{U}]$ consist of elements $p$ of the form
$p=\langle t_1,...,t_n,\langle\kappa,B\rangle\rangle$.
 For every  $1\leq i\leq n $, $t_i$ is either an ordinal
 $\kappa_i$ if $ o^{\vec{U}}(\kappa_i)=0$
 or a pair $\langle\kappa_i,B_i\rangle$  if \ $o^{\vec{U}}(\kappa_i)>0$.
\begin{enumerate}
\item $B\in\cap\vec{U}(\kappa)$, \  $\min(B)>\kappa_n$.
    \item  For every  $1\leq i\leq n$.
    \begin{enumerate}
    \item $\langle\kappa_1,...,\kappa_n\rangle\in [\kappa]^{<\omega}$ (increasing finite sequence below $\kappa$).
    \item $B_i\in \cap\vec{U}(\kappa_i)$.
    \item  $\min(B_i)>\kappa_{i-1}$ \ $(i>1)$.
    \end{enumerate}
\end{enumerate}

\end{definition}
\begin{definition}\label{Magidor-order}
 For $p=\langle t_1,t_2,...,t_n,\langle\kappa,B\rangle\rangle,q=\langle s_1,...,s_m,\langle\kappa,C\rangle\rangle\in \Mfor$ , define  $p \leq q$ ($q$ extends $p$) iff:
\begin{enumerate}
    \item $n \leq m$.
    \item $B \supseteq C$.
    \item $\exists 1 \leq i_1 <...<i_n \leq m$ such that for every $1 \leq j \leq m$:
    \begin{enumerate}
        \item If $\exists 1\leq r\leq n$ such that $i_r=j$ then $\kappa(t_r)=\kappa( s_{i_r})$ and $C(s_{i_r})\subseteq B(t_r)$.
        \item Otherwise $\exists \ 1 \leq r \leq n+1$ such that $ i_{r-1}<j<i_{r}$ then 
        \begin{enumerate}
        \item $\kappa(s_j) \in B(t_r)$.
        \item $B(s_j)\subseteq B(t_r)\cap \kappa(s_j)$.
        \item $o^{\vec{U}}(s_j)<o^{\vec{U}}(t_r)$.
        \end{enumerate}
    \end{enumerate}
\end{enumerate}
We also use``p directly extends q", $p \leq^{*} q$ if:
\begin{enumerate}
    \item $p \leq q$
    \item $n=m$
\end{enumerate}
\end{definition}
Let us add some notation, for a pair $t=\langle \alpha, X\rangle$ we denote by $\kappa(t)=\alpha,\ B(t)=X$. If $t=\alpha$ is an ordinal then $\kappa(t)=\alpha$ and $B(t)=\emptyset$.

For a condition $p=\langle t_1,...,t_n,\langle \kappa,B\rangle\rangle\in\Mfor$ we denote
$n=l(p)$, $p_i=t_i$, $B_i(p)=B(t_i)$ and $\kappa_i(p)=\kappa(t_i)$ for any $1\leq i\leq l(p)$, $t_{l(p)+1}=\langle\kappa,B\rangle$, $t_0=0$. Also denote
$$\kappa(p)=\{\kappa_i(p)\mid i\leq l(p)\}\text{ and }B(p)=\bigcup_{i\leq l(p)+1}B_i(p)$$

\begin{remark}
Condition 3.b.iii is not essential, since the set
$$\Big\{p\in\Mfor\mid \forall i\leq l(p)+1. \forall\alpha\in B_{i}(p).o^{\vec{U}}(\alpha)<o^{\vec{U}}(\kappa_i(p))\Big\}$$
is a dense subset of $\Mfor$ and the order between any two elements of this dense subsets automatically satisfy 3.b.iii.
\end{remark}

\begin{definition}\label{end extension}
  Let $p\in\Mfor$. For every $ i\leq l(p)+1$, and $\alpha\in B_{i}(p)$ with $o^{\vec{U}}(\alpha)>0$, define 
$$p^{\frown}\langle\alpha\rangle=\langle p_1,...,p_{i-1},\langle\alpha,B_{i}(p)\cap\alpha\rangle,\langle\kappa_{i}(p),B_{i}(p)\setminus(\alpha+1)\rangle,p_{i+1},...,p_{l(p)+1}\rangle$$
If $o^{\vec{U}}(\alpha)=0$, define
$$p^{\frown}\langle\alpha\rangle=\langle p_1,...,p_{i-1},\alpha,\langle\kappa_{i}(p),B_{i}(p)\setminus(\alpha+1)\rangle,...,p_{l(p)+1}\rangle$$
For $\langle\alpha_1,...,\alpha_n\rangle\in[\kappa]^{<\omega}$ define recursively,
$$p^{\frown}\langle\alpha_1,...,\alpha_n\rangle=(p^{\frown}\langle\alpha_1,...,\alpha_{n-1}\rangle)^{\frown}\langle\alpha_n\rangle$$

\end{definition}
\begin{proposition}
Let $p\in\Mfor$. If $p^{\frown}\vec{\alpha}\in\Mfor$, then it is the minimal extension of $p$ with stem $$\kappa(p)\cup\{\vec{\alpha}_1,...,\vec{\alpha}_{|\vec{\alpha}|}\}$$ Moreover, $p^{\frown}\vec{\alpha}\in\Mfor$ iff for every $i\leq |\vec{\alpha}|$ there is $j\leq l(p)$ such that:
\begin{enumerate}
   \item $\vec{\alpha}_i\in(\kappa_j(p),\kappa_{j+1}(p))$.
\item  $o^{\vec{U}}(\vec{\alpha}_i)<o^{\vec{U}}(\kappa_{j+1})$.
\item $B_{j+1}(p)\cap\vec{\alpha}_i\in\cap\vec{U}(\vec{\alpha}_i)$.

\end{enumerate}
\end{proposition}
$\blacksquare$
\vskip 0.2 cm
Note that if we add a pair of the form $\langle \alpha , B\cap\alpha\rangle$ then in $B\cap\alpha$ there might be many ordinals which are irrelevant to the forcing. Namely, ordinals $\beta\in B\cap\alpha$ with $o^{\vec{U}}(\beta)\geq o^{\vec{U}}(\alpha)$, such ordinals cannot be added to the sequence.

\begin{definition}
Let $p\in \Mfor$, define for every $i\leq l(p)$
 $$p\restriction\kappa_i(p)= \langle p_1,...,p_i\rangle \ and \ p\restriction(\kappa_i(p),\kappa)=\langle p_{i+1},...,p_{l(p)+1}\rangle$$
 Also, for $\lambda$ with $o^{\vec{U}}(\lambda)>0$ define
 $$\Mfor\restriction\lambda=\{p\restriction\lambda\mid p\in\Mfor \ and \ \lambda \ appears \ in \ p\}$$
 $$\Mfor\restriction(\lambda,\kappa)=\{p\restriction (\lambda,\kappa)\mid p\in\Mfor \ and \ \lambda \ appears \ in \ p\}$$
\end{definition}
Note that $\Mfor\restriction\lambda$ is just Magidor forcing on $\lambda$ and $\Mfor\restriction (\lambda,\kappa)$ is a subset of $\Mfor$. The following decomposition is straightforward.
\begin{proposition}\label{dec} Let $p\in\Mfor$ and $\langle\lambda,B\rangle$ a pair in $p$. Then
$$\Mfor/p\simeq \Big(\Mfor\restriction \lambda\Big)/\Big(p\restriction\lambda\Big)\times\Big(\Mfor\restriction(\lambda,\kappa)\Big)/\Big(p\restriction(\lambda,\kappa)\Big)$$
\end{proposition}
\begin{remark}\label{lift}
When considering $\vec{U}$ in some model $V\subseteq N\subseteq V[C_G\cap\lambda]$, since we added generic sequences, not all of the measures in $\vec{U}$ remain measures in $N$. However, each measure $U(\xi,i)$ for $\lambda<\xi\leq\kappa$ and $i<o^{\vec{U}}(\xi)$ generates a normal measure $W(\xi,i)$ over $\xi$ such that $\vec{W}=\l W(\xi,i)\mid \lambda<\xi\leq\kappa, \ i<o^{\vec{U}}(\xi)\r$ is a coherent sequence. Since  $\Mfor\restriction(\lambda,\kappa)$ is a dense subset of $\mathbb{M}[\vec{W}]$, forcing over $N$ with  $\Mfor\restriction(\lambda,\kappa)$ is the same as forcing with $\mathbb{M}[\vec{W}]$. 
\end{remark}
\begin{proposition}
Let $p\in\Mfor$ and $\langle\lambda,B\rangle$ a pair in $p$. Then the order $\leq^*$ in the forcing $\Big(\Mfor\restriction(\lambda,\kappa)\Big)/\Big(p\restriction(\lambda,\kappa)\Big)$ is $\delta$-directed where $\delta=\min\{\nu>\lambda\mid o^{\vec{U}}(\nu)>0\}$. Meaning that for every $X\subseteq \Mfor\restriction (\lambda,\kappa)$ such that $|X|<\delta$ and for every $q\in X, \ p\leq^* q$, there is an  $\leq^*$-upper bound for $X$.
\end{proposition}
 \begin{lemma}
$\Mfor$ satisfies $\kappa^+$-c.c.
\end{lemma}
The following is known as the Prikry condition:
\begin{lemma}\label{prikrycondition}
$\Mfor$ satisfy the Prikry condition i.e. for any statement in the forcing language $\sigma$ and any $p\in\Mfor$ there is $p\leq^*p^*$ such that $p^*||\sigma$ i.e. either $p^*\Vdash\sigma$ or $p\Vdash\neg\sigma$.
\end{lemma} 

The next lemma can be found in \cite{ChangeCofinality}:
\begin{lemma}\label{MagLemma}
  Let $G\subseteq \Mfor$ be generic and
suppose that $A\in V[G]$ is such that $A\subseteq V_\alpha$. Let $p\in G$ and $\l\lambda,B\r$ a pair in $p$ such that $\alpha<\lambda$, then $A\in V[G\restriction\lambda]$.
\end{lemma}
\pr  Consider the decomposition \ref{dec} $p=\langle q,r\r$, where $q\in \Mfor\restriction \lambda$ and $r\in\Mfor\restriction (\lambda,\kappa)$
Work in $V[G\restriction\lambda]$, 
Let $\lusim{A}$ be a $\Mfor\restriction (\lambda,\kappa)$-name for $A$. For every $x\in V_\alpha$ use 
the Prikry condition \ref{prikrycondition}, to find $r\leq^* r_x$ such that $r_x$ decides the statement $r\in \lusim{A}$.
By definition of $\lambda$ and proposition \ref{decprop}, the forcing $\Mfor\restriction(\lambda,\kappa)$ is $|V_\alpha|^+$-directed with the $\leq^*$ order.
Hence there is $r\leq^* r^*$ such that $p_x\leq^* p^*$ for every $x\in V_{\alpha}$. By density, we can find such $r^*\in G\restriction(\lambda,\kappa)$. It follows that $A=\{x\in V_\alpha\mid r^*\Vdash x\in \lusim{A}\}$ is definable in $V[G\restriction \lambda]$.$\blacksquare$
\begin{corollary}
$\Mfor$ preserves all cardinals.
\end{corollary}
 \begin{definition}
  Let $G\subseteq \Mfor$ be generic, define the \textit{Magidor club}
  $$C_{G}=\{ \nu \mid \exists p\in G \exists i\leq l(p)\text{ s.t. } \nu=\kappa_i(p)\}$$

\end{definition}
  We will abuse notation by sometimes considering $C_G$ as a the canonical enumeration of the set $C_G$. The set $C_{G}$ is closed and unbounded in $\kappa$, therefore, the order type of $C_{G}$ determines the cofinality of $\kappa$ in $V[G]$. The next propositions can be found in \cite{Gitik2010}.
 \begin{proposition}\label{decprop}
 Let $G\subseteq\Mfor$ be generic. Then $G$ can be reconstructed from $C_{G}$ as follows
     $$ G=\{p\in\Mfor\mid (\kappa(p)\subseteq C_{G}) \wedge (C_{G}\setminus\kappa(p)\subseteq B(p))\}$$
      In particular $V[G]=V[C_{G}]$.
\end{proposition}
\begin{proposition}\label{genericproperties}
Let $G\subseteq\Mfor$ be generic. 
\begin{enumerate}
    \item $C_G$ is a club at $\kappa$.
    \item For every $\delta\in C_G$, $o^{\vec{U}}(\delta)>0$ iff $\delta\in Lim(C_G)$.
    \item For every $\delta\in Lim(C_G)$, and every $A\in \cap\vec{U}(\delta)$, there is $\xi<\delta$ such that $C_G\cap(\xi,\delta)\subseteq A$.
    
    \item If $\l \delta_i\mid i<\theta\r$ is an increasing sequence of elements of $C_G$, let $\delta^*=\sup_{i<\theta}\delta_i$, then $o^{\vec{U}}(\delta^*)\geq\limsup_{i<\theta}o^{\vec{U}}(\delta_i)+1$.\footnote{ For a sequence of ordinals $\l \rho_j\mid j<\gamma\r$, $\limsup_{j<\gamma}\rho_j=\min\{\sup_{i<j<\gamma}\rho_j\mid i<\gamma\}$.}
    \item Let $\delta\in Lim(C_G)$ and let $A$ be a positive set, $A\in (\cap\vec{U}(\delta))^+$. i.e. $
\delta\setminus A\notin \cap\vec{U}(\delta)$. \footnote{Equivalently, if there is some $i<o^{\vec{U}}(\delta)$ such that $A\in U(\delta,i)$.} Then, $\sup(A\cap C_G)=\delta$.
\item If $A\subseteq V_\alpha$, then $A\in V[C_G\cap\lambda]$, where $\lambda=\max(Lim(C_G)\cap\alpha+1)$.
\item For every $V$-regular cardinal $\alpha$, if $cf^{V[G]}(\alpha)<\alpha$ then $\alpha\in Lim(C_G)$.
\end{enumerate}

\end{proposition}
\pr 
$(1),(2),(3)$ can be found in \cite{Gitik2010}.

To see $(4)$, use closure of $C_G$, and find $q\in G$ such that $\delta^*$ appears in $q$. Since there are only finitely many ordinals in $q$, there is some $i<\theta$ such that for every $j>i$, $\delta_j$ does not appear in $q$. By \ref{Magidor-order}, since every such $\delta_j$ appear in some $q_j\in G$ which is compatible with $q$, $o^{\vec{U}}(\delta_j)<o^{\vec{U}}(\delta^*)$.
Hence $$\limsup_{j<\theta}o^{\vec{U}}(\delta_j)+1\leq \sup_{i<j<\theta}o^{\vec{U}}(\delta_j)+1\leq o^{\vec{U}}(\delta^*)$$
For $(5)$, let $\rho<\delta$. Each condition $p$, such that $\delta=\kappa_i(p)$ for some $i\leq l(p)+1$, must satisfy that $\sup(A\cap B_i(p))=\delta$. Hence we can extend $p$ using an element of $A\cap B_i(p)$ above $\rho$. By density, $\sup(A\cap C_G)\geq\rho$. Since $\rho$ is general, $\sup(A\cap C_G)=\delta$.

$(6)$ is a direct corollary of \ref{MagLemma}. As for $(7)$, assume that $cf^{V[G]}(\alpha)<\alpha$, and let $X\subseteq\alpha$ be a club such that $\otp(X)=cf^{V[G]}(\alpha)$. Then $X\in V[G]\setminus V$. Let $\lambda=\max(Lim(C_G)\cap\alpha+1)$, then $\lambda\leq\alpha$. By $(6)$, $X\in V[C_G\cap \lambda]$. Toward a contradiction, assume that $\lambda<\alpha$, then the forcing $\Mfor\restriction \lambda$ is $\alpha$-c.c., but $cf^{V[C_G\cap\lambda]}(\alpha)<\alpha$, contradiction.
$\blacksquare$

The Mathias-like criteria for Magidor forcing is due to Mitchell \cite{MitchellMathias}:
\begin{theorem}\label{MathiasCriteria}
Let $U$ be a coherent sequence and assume that $c:\alpha\rightarrow \kappa$ is an increasing function. Then $c$ is $\Mfor$ generic iff:
\begin{enumerate}
 \item $c$ is continuous.
 \item $c\restriction\beta$ is $\mathbb{M}[\vec{U}\restriction\beta]$ generic for every $\beta\in \Lim(\alpha)$.
 \item $X\in \cap\vec{U}(\kappa)$ iff $\exists \beta<\kappa \ Im(c)\setminus\beta\subseteq X$.
\end{enumerate}
\end{theorem}
An equivalent formulation of the Mathias criteria is to require that for every $\beta\in \Lim(\alpha)$, and for every $X\in \cap\vec{U}(c(\beta))$, there is $\xi<\beta$ such that $c''(\xi,\beta)\subseteq X$.

For an additional proof of \ref{MathiasCriteria}, we refer the reader to the last section, where we define a forcing notion $\mathbb{M}_f[\vec{U}]$, which generalizes $\Mfor$, and prove in \ref{MathisForMfU} a Mathias-like criteria for it. 

\begin{proposition}\label{indc}
 Let $G\subseteq \Mfor$ be $V$-generic filter and $C_{G}$ the corresponding Magidor sequence. Let $p\in G$, then for every $i\leq l(p)+1$ 
 \begin{enumerate}
     \item If  $o^{\vec{U}}(\kappa_i(p))\leq \kappa_i(p)$,
    $$\otp( [\kappa_{i-1}(p),\kappa_i(p))\cap C_{G} )=\omega^{o^{\vec{U}}(\kappa_i(p))}$$
    \item If $o^{\vec{U}}(\kappa_i(p))\geq \kappa_i(p)$, then 
    $$\otp( [\kappa_{i-1}(p),\kappa_i(p))\cap C_{G} )=\kappa_i(p)$$
    
 \end{enumerate}
 \end{proposition}
 \pr We prove $(1)$ by induction on $\kappa_i(p)$.
 If $\kappa_i(p)=0$, then $C_G\cap [\kappa_{i-1}(p),\kappa_i(p))=\{\kappa_{i-1}(p)\}$. Hence $$\otp(C_G\cap [\kappa_{i-1}(p),\kappa_i(p)))=1=\omega^0=\omega^{o^{\vec{U}}(\kappa_i(p))}$$
 
Assume the lemma holds for any $\delta<\kappa_i(p)$. If $o^{\vec{U}}(\kappa_i(p))=\alpha+1\leq\kappa_i(p)$, then the set $X=\{\beta<\kappa_i(p)\mid o^{\vec{U}}(\beta)=\alpha\}\in U(\kappa_i(p),\alpha)$, hence by proposition \ref{genericproperties}, $$\sup (X\cap C_G\cap [\kappa_{i-1}(p),\kappa_i(p)))=\kappa_i(p)$$  We claim that $\otp(X\cap C_G\cap [\kappa_{i-1}(p),\kappa_i(p))=\omega$. Otherwise, let $\rho<\kappa_i(p)$ be such that $\rho$ is a limit point of $X\cap C_G\cap [\kappa_{i-1}(p),\kappa_i(p))$. Again by proposition \ref{genericproperties}, $$o^{\vec{U}}(\rho)\geq \limsup(o^{\vec{U}}(\xi)\mid \xi\in X\cap C_G\cap [\kappa_{i-1}(p),\kappa_i(p)))=\alpha+1$$ Contradicting \ref{Magidor-order}. Let $\l \delta_n\mid n<\omega\r$ be the increasing enumeration of $X\cap C_G\cap [\kappa_{i-1}(p),\kappa_i(p))$. By induction hypothesis, for every $n<\omega$, $\otp(C_G\cap[\delta_n,\delta_{n+1}))=\omega^\alpha$. Hence, $$\otp( C_G\cap [\kappa_{i-1}(p),\kappa_i(p))=\omega^{\alpha+1}$$
For limit $o^{\vec{U}}(\kappa_i(p))$, use proposition \ref{genericproperties}(5), to see that the sequence $\l\delta_\alpha\mid \alpha<o^{\vec{U}}(\kappa_i(p))\r$ where $$\delta_\alpha=\min\{\rho\in C_G\cap [\kappa_{i-1}(p),\kappa_i(p))\mid  o^{\vec{U}}(\rho)=\alpha\}$$
is well defined. $x=\sup(\delta_\alpha\mid \alpha<\theta)\leq \kappa_i(p)$ is an element of $C_G$, and by proposition \ref{genericproperties}(4), $o^{\vec{U}}(x)\geq o^{\vec{U}}(\kappa_i(p))$, hence $x=\kappa_i(p)$. For every $\alpha<o^{\vec{U}}(\kappa_i(p))$, $\otp(C_G\cap [\kappa_i(p),\delta_\alpha))=\omega^\alpha$, since $p^{\smallfrown}\l \delta_\alpha\r\in G$ and  by induction hypothesis. It follows that $$\otp(C_G\cap [\kappa_{i-1}(p),\kappa_i(p))=\sup_{\alpha<o^{\vec{U}}(\kappa_i(p))}(\otp(C_G\cap[\kappa_{i-1}(p),\delta_\alpha))=\sup_{\alpha<o^{\vec{U}}(\kappa_i(p))}\omega^\alpha=\omega^{o^{\vec{U}}(\kappa_i(p))}$$

For $(2)$, use $(1)$, and the limit stage to conclude that if $o^{\vec{U}}(\kappa_i(p))=\kappa_i(p)$, then $$\otp(C_G\cap [\kappa_{i-1}(p),\kappa_i(p))=\kappa_i(p)$$
If $o^{\vec{U}}(\kappa_i(p))>\kappa_i(p)$, then $\{\alpha<\kappa_i(p))\mid o^{\vec{U}}(\alpha)=\alpha\}\in U(\kappa_i(p),\kappa_i(p))$, hence by proposition \ref{genericproperties}, there are unboundedly many $\alpha\in C_G\cap[\kappa_{i-1}(p),\kappa_i(p))=:Y$ such that $o^{\vec{U}}(\alpha)=\alpha$. 
Hence $$\kappa_i(p)=\sup(Y)=\sup(\otp( C_G\cap[\kappa_{i-1}(p),\alpha)\mid \alpha\in Y)\leq \kappa_i(p)$$
So equality holds.$\blacksquare$

  Proposition \ref{indc} suggests a connection between the index in $C_G$ of ordinals appearing in $p$ 
  and Cantor normal form.

  \begin{definition}\label{DefLimitOrder}
 Let $p\in G$. For each $i\leq l(p)$ define
   $$\gamma_i(p)=\sum_{j=1}^{i}\omega^{o^{\vec{U}}(\kappa_j(p))}$$ 
Also for an ordinal $\alpha$, denote by $o_L(\alpha)=\gamma_n$ where $\alpha=\sum_{i=1}^n\omega^{\gamma_i}\cdot m_i$ is the Cantor normal form and $\gamma_1>\gamma_2...>\gamma_n$.
   \end{definition}
 \begin{corollary}\label{IndCG}
 Let $G\subseteq\Mfor$ be $V$-generic and $C_G$ the corresponding Magidor sequence.
 \begin{enumerate}
     \item If $p\in G$, then for every $1\leq i\leq l(p)$
 $$p\Vdash \underaccent{\sim}{C}_G(\gamma_i(p))=\kappa_i(p)$$
     \item For every $\alpha<\otp(C_G)$, $o^{\vec{U}}(C_G(\alpha))=o_L(\alpha)$
 \end{enumerate}
 \end{corollary}
\pr This is directly from \ref{indc}.$\blacksquare$

For more details and basic properties of Magidor forcing see \cite{ChangeCofinality},\cite{Gitik2010} or \cite{TomMoti}. 

We are going to handle subsequences of the generic club, the following simple definition will turn out being useful.
\begin{definition}\label{Index}
 Let $X,X'$ be sets of ordinals such that $X'\subseteq X\subseteq On$. Let $\alpha=otp(X,\in)$ be the order type of $X$ and $\phi:\alpha\rightarrow X$ be the order isomorphism witnessing it. The indices of $X'$ in $X$ are $$I(X',X)=\phi^{-1''}X'=\{\beta<\alpha\mid \phi(\beta)\in X'\}$$
\end{definition}

In the last part of the proof we will need the definition of quotient forcing.
\begin{definition}
 Let $\underaccent{\sim}{C'}$ be a $\Mfor$-name for a subset of $C_G$, and let $C'\subseteq C_G$ such that $\underaccent{\sim}{C'}_G=C'$. Define $\mathbb{P}_{\underaccent{\sim}{C'}}$, the complete subalgebra of $\l RO(\Mfor),\leq_B\r$\footnote{ $RO(\Mfor)$ is the set of all regular open cuts of $\Mfor$(see for example \cite[Thm. 14.10]{Jech2003}), as usual we identify $\Mfor$ as a dense subset of $RO(\Mfor)$. The order $\leq_B$ is in the standard position of Boolean algebras orders i.e. $p\leq_B q$ means $p\Vdash q\in \hat{G}$.} generated by the conditions $X=\{||\alpha\in \underaccent{\sim}{C'}||\mid \alpha<\kappa\}$. is\end{definition}
By \cite[15.42]{Jech2003}, $V[C']=V[H]$ for some $V$-generic filter $H$ of $\mathbb{P}_{\underaccent{\sim}{C'}}$. In fact $$C'=\{\alpha<\kappa\mid ||\alpha\in \underaccent{\sim}{C'}||\in X\cap H\}$$
\begin{definition}
Define the function $\pi:\Mfor\rightarrow \mathbb{P}_{\underaccent{\sim}{C'}}$ by 
$$\pi(p)=\inf(b\in \mathbb{P}_{\underaccent{\sim}{C'}}\mid p\leq_B b)$$
\end{definition}
It not hard to check that $\pi$ is a projection i.e. 
\begin{enumerate}
    \item $\pi$ is order preserving.
    \item $\forall p\in\Mfor.\forall q\leq_B\pi(p).\exists p'\geq p.\pi(p')\leq_B q$.
    \item $Im(\pi)$ is dense in $\mathbb{P}_{\underaccent{\sim}{C'}}$.
    \end{enumerate}
\begin{definition}
Let $\pi:\mathbb{P}\rightarrow\mathbb{Q}$ be any projection, let $H\subseteq\mathbb{Q}$ be $V$-generic, define
$$\mathbb{P}/H=\pi^{-1''}H$$
\end{definition}
We abuse notation by defining $\Mfor/C'=\Mfor/H$, where $H$ is some generic for $\mathbb{P}_{\underaccent{\sim}{C'}}$ such that $V[H]=V[C']$. It is known that if $G$ is $V[C']$-generic for $\Mfor/C'$ then $G$ is $V$ generic for $\Mfor$ and $\bar{\pi''G}=H$, hence $V[G]=V[C'][G]$.

\section{Magidor forcing with $o^{\vec{U}}(\kappa)\leq\kappa$}

Assume that $o^{\vec{U}}(\kappa)\leq\kappa$.  Let $G\subseteq\Mfor$ be a $V$-generic filter, and let $p\in G$. By proposition \ref{indc}, $otp(C_G\cap (\kappa_{l(p)}(p),\kappa))=\omega^{o^{\vec{U}}(\kappa)}$. Hence, $$(1.1) \ \ \ cf^{V[G]}(\kappa)=cf^{V[G]}(\omega^{o^{\vec{U}}(\kappa)})$$
\begin{corollary}
\begin{enumerate}
    \item If $o^{\vec{U}}(\kappa)<\kappa$, then $\kappa$ is singular in $V[G]$.
    \item If $o^{\vec{U}}(\kappa)=\kappa$, then $cf^{V[G]}(\kappa)=\omega$.
\end{enumerate}

\end{corollary}
\pr $(1)$ follows directly from equation $(1.1)$. For $(2)$,
the set $E=\{\alpha<\kappa\mid o^{\vec{U}}(\alpha)<\alpha\}\in\cap\vec{U}(\kappa)$. Hence, by proposition \ref{genericproperties} find $\rho<\kappa$ such that $C_G\setminus \rho\subseteq E$.
In $V[G]$ consider the sequence: $\alpha_0=\min(C_G\setminus\rho)$, then $\alpha_{n+1}=C_G(\alpha_n)$. This is a well defined sequence of ordinals below $\kappa$ since $\otp(C_G)=\kappa$. Also, since $\{\alpha<\kappa\mid \omega^\alpha=\alpha\}\in\cap\vec{U}(\kappa)$, there is $n<\omega$, such that for every $m\geq n$, $o^{\vec{U}}(\alpha_{m+1})=\alpha_m$.

To see that $\alpha^*:=\sup_{n<\omega}\alpha_n=\kappa$, assume otherwise, then by closure of $C_G$, $\alpha^*\in C_G$. Also $\alpha^*>\rho$, hence $o^{\vec{U}}(\alpha^*)<\alpha^*$. By proposition \ref{genericproperties}(4), $$o^{\vec{U}}(\alpha^*)\geq \limsup_{n<\omega}o^{\vec{U}}(\alpha_n)+1=\sup_{n<\omega}\alpha_n=\alpha^*$$
contradiction.$\blacksquare$

 If $o^{\vec{U}}(\kappa)\leq\kappa$. We can decompose every set $A\in\cap\vec{U}(\kappa)$ in a very canonical way:
\begin{proposition}\label{simple case} 
Assume that $o^{\vec{U}}(\kappa)\leq\kappa$.  Let $A\in
\cap\vec{U}(\kappa)$.
\begin{enumerate}
\item For every $i<\kappa$ define $A_i=\{\nu\in A\mid o^{\vec{U}}(\nu)=i\}$. Then $A=\underset{i<\kappa}{\biguplus}A_i$ and $A_i\in U(\kappa,i)$. 
 \item There exists $A^*\subseteq A$ such that:
 \begin{enumerate}
  \item $A^*\in\cap\vec{U}(\kappa)$
  \item For every $0<j<o^{\vec{U}}(\kappa)$ and $\alpha\in A_j^*$,  $A^*\cap\alpha\in\cap\vec{U}(\alpha)$.
 \end{enumerate}
\end{enumerate}
\end{proposition}
\pr
1. Note that $X_i:=\{\nu<\kappa\mid o^{\vec{U}}(\nu)=i\}\in U(\kappa,i)$ and $A_i=X_i\cap A\in U(\kappa,i)$. Moreover, every $\alpha<\kappa$, $o^{\vec{U}}(\alpha)<\kappa$, since there are at most $2^{2^{\alpha}}<\kappa$ measures over $\alpha$.\\ 
2. For any $i<o^{\vec{U}}(\kappa)$,
$$Ult(V,U(\kappa,j))\models A=j_{U(\kappa,j)}(A)\cap\kappa\in \underset{i<j}{\bigcap}U(\kappa,i)$$
Coherency of the sequence implies that $A':=\{\alpha<\kappa\mid A\cap\alpha\in\cap\vec{U}(\alpha)\}\in U(\kappa,j)$, this is for every $j<o^{\vec{U}}(\kappa)$. \\
 Define inductively $A^{(0)}=A$, $A^{(n+1)}=A^{'(n)}$. By definition, $\forall\alpha\in A^{(n+1)}_j$,  $A^{(n)}\cap\alpha\in\cap\vec{U}(\alpha)$. Define $A^*=\underset{n<\omega}{\bigcap}A^{(n)}\in\cap\vec{U}(\kappa)$, this set
has the required property.
$\blacksquare$

\subsection{Extension Types}

By convention, for a set of ordinals $B$, $[B]^{<\alpha}$ is the set of increasing sequences of length less than $\alpha$ of ordinals in $B$, $[B]^{[
<\alpha]}$ is the set of not necessarily increasing sequences of length less than $\alpha$ of ordinals in $B$. For sets of ordinals $B_i$ for $1\leq i\leq n$, let $\prod_{i=1}^nB_i$ be the set of increasing sequence $\l\alpha_1,..,\alpha_n\r$ such that $\alpha_i\in B_i$. For double indexed sets $B_{i,j}$ for $1\leq i\leq n, \ 1\leq j \leq m$, the set $\prod_{i=1}^n\prod_{j=1}^nB_{i,j}$ is viewed as a product of single indexed sets using the left lexicographical order.
\begin{definition}
Let $p\in \Mfor$. Define
\begin{enumerate}
\item For every $i\leq l(p)+1$, let $B_{i,\alpha}(p)=B_i(p)\cap X_\alpha$, where $X_\alpha:=\{\beta<\kappa\mid o^{\vec{U}}(\beta)=\alpha\}$ are the sets defined in \ref{simple case}.
 
 \item $Ex(p)=\prod^{l(p)+1}_{i=1}[o^{\vec{U}}(\kappa_i(p))]^{[<\omega]}$. 
 \item If $X\in Ex(p)$, then $X$ is of the form $\langle X_1,...,X_{n+1}\rangle$. Denote $x_{i,j}$, the $j$-th element of $X_i$, for $1\leq j\leq |X_i|$ and $mc(X)$ is the last element of $X$ and $l(X)=\sum_{i=1}^{n+1}|X_i|$.
\item Let $X\in Ex(p)$, then
$$\vec{\alpha}=\langle\vec{\alpha_1},...,\vec{\alpha_{l(p)+1}}\rangle \in  \overset{l(p)+1}{\underset{ i=1}{\prod}}\overset{|X_i|}{\underset{ j=1}{\prod}} B_{i , x_{i,j}}(p)=:X(p)$$
 call $X$ an \textit{extension-type of $p$} and $\vec{\alpha}$ is of \textit{type $X$}, note that $\vec{\alpha}$ is an increasing sequence of ordinals.
 \end{enumerate}
 \end{definition}
 The idea of extension types is simply to classify extensions of $p$ according to the measures from which the ordinals added to the stem of $p$ are chosen.
Note that if $o^{\vec{U}}(\kappa)=\lambda<\kappa$ then there is a bound on the number of extension types, $|Ex(p)|<\min\{\nu>\lambda\mid o^{\vec{U}}(\nu)>0\}$.

 By proposition \ref{simple case} any $p\in\Mfor$ can be extended to $p\leq^*p^*$ such that for every $X\in Ex(p)$ and any $\vec{\alpha}\in X(p)$, $p^{\frown}\vec{\alpha}\in \Mfor$. Let us move to this dense subset of $\Mfor$. 
\begin{proposition}\label{Max}
 Let $p\in \Mfor$ be any condition and
 $p\leq q\in \Mfor$. Then there exists unique $X\in Ex(p)$ and $\vec{\alpha}\in X(p)$ such that $p^{\frown}\vec{\alpha}\leq^* q$. Moreover, for every $X\in Ex(p)$ the set $\{p^{\frown}\vec{\alpha}\mid\vec{\alpha}\in X(p)\}$ form a maximal antichain above $p$. 
\end{proposition}
\pr The first part is trivial. We will prove that $\{p^{\smallfrown}\vec{\alpha}\mid\vec{\alpha}\in X(p)\}$ form a antichain above $p$, by induction on $l(X)$. For $l(X)=1$, we merely have some $X(p)=B_{i,\xi}(p)\in U(\kappa_i(p),\xi)$.
To see it is an antichain, let $\beta_1<\beta_2$ are in $X(p)$. Toward a contradiction, assume that $p^{\smallfrown}\beta_1,p^{\smallfrown}\beta_2\leq q$, then $\beta_1$ appears in a pair in $q$ and is added between $\kappa_{i-1}(p)$ and $\beta_2$, so by definition \ref{Magidor-order}, it must be that $\xi=o^{\vec{U}}(\beta_1)<o^{\vec{U}}(\beta_2)=\xi$ contradiction.

To see it is maximal, fix $q\geq p$ and let $\vec{\alpha}$ be such that $p^{\frown}\vec{\alpha}\leq^*q$. Consider the type of $\vec{\alpha}$, $$Y\in Ex(p)$$, then $\vec{\alpha}\in Y(p)$. In $Y_i$ let $j$ be the minimal such that $y_{i,j}\geq\xi$.
If $y_{i,j}=\xi$ then $q\geq p^{\frown}\langle\alpha_{i,j}\rangle\in X(p)$ and we are done. Otherwise, $y_{i,j}>\xi$, then one of the pairs in $q$ is of the form
$\langle\alpha_{i,j},B\rangle$ where $B\in \cap\vec{U}(\alpha_{i,j})$ and $B\subseteq B_i(p)$. Any $\alpha\in B\cap B_{i,\xi}(p)$, will satisfy that $p^{\frown}
\l\alpha\r\in X(p)$ and $p^{\frown}\l\alpha\r,q\leq q^{\frown}\l\alpha\r$.

Assume that the claim holds for $l(X)=n$, and let $X\in Ex(p)$ be such that $l(X)=n+1$. Let $\vec{\alpha},\vec{\beta}\in X(p)$ be distinct, if for some $x_{i,j}\neq mc(X)$ we have $\alpha_{i,j}\neq\beta_{i,j}$ apply the induction to $X\setminus mc(X)$ to see that $p^{\frown}\vec{\alpha}\setminus \alpha^*,p^{\frown}\vec{\beta}\setminus\beta^*$ are incompatible, hence $p^{\frown}\vec{\alpha},p^{\frown}\vec{\beta}$ are incompatible. If $\vec{\alpha}\setminus \alpha^*=\vec{\beta}\setminus\beta^*$, then $\alpha^*\neq\beta^*$ and by the case $n=1$ we are done. To see it is maximal, let $q\geq p$ apply the induction to $X'$ which is the extension type obtained from $X$ by removing $mc(X)$ to find $\vec{\alpha}\in X'(p)$ such that $p^{\frown}\vec{\alpha}$ is compatible with $q$ and let $q'$ be a common extension. Again by the case $n=1$, there is $\l\alpha\r\in mc(X)(p^{\frown}\vec{\alpha})$ such that $p^{\frown}\vec{\alpha}^{\frown}\l\alpha\r$ and $q'$ are compatible. 
 $\blacksquare$
\begin{definition}
  Let $U_1,...,U_n$ be ultrafilters on a $\kappa_1\leq...\leq\kappa_n$ respectively, define recursively the ultrafilter $\prod_{i=1}^nU_i$ over $ \prod_{i=1}^n\kappa_i$, as follows: for $B\subseteq \prod_{i=1}^n\kappa_i$
  $$B\in \prod_{i=1}^nU_i\leftrightarrow \{\alpha_1<\kappa_1\mid B_{\alpha_1}\in \prod_{i=2}^nU_i\}\in U_1$$
  where $B_\alpha=B\cap\Big(\{\alpha\}\times \prod_{i=2}^n\kappa_i\Big)$.
\end{definition}
\begin{proposition}\label{Normgenerate}
If $U_1,...,U_n$ are normal ultrafilter, then $\prod_{i=1}^nU_i$ is generated by sets of the form $A_1\times...\times A_n$ such that $A_i\in U_i$.
\end{proposition}
\pr By induction of $n$, for $n=1$ there is nothing to prove. Assume that the proposition holds for $n-1$, and let $B\in \prod_{i=1}^nU_i$. By definition,  $A_1=\{\alpha_1<\kappa_1\mid B_{\alpha_1}\in\prod_{i=2}^n U_i\}\in U_1$, and by the induction hypothesis each $B_{\alpha_1}$ contains a set of the form $A_{2,\alpha_1}\times...\times A_{n,\alpha_1}$. By normality,  $A_i:=\Delta_{\alpha\in A_1}A_{i,\alpha}\in U_i$. Consider $\l\alpha_1,...,\alpha_n\r\in A_1\times....\times A_n$, by convention, for each $2\leq i\leq n$, $\alpha_1\leq \alpha_i$, and by definition of diagonal intersection, $\alpha_i\in A_{i,\alpha_1}$, hence $\l\alpha_2,...,\alpha_n\r\in A_{2,\alpha_1}\times...\times A_{n,\alpha_1}\subseteq B_{\alpha_1}$. It follows by the definition of $B_{\alpha_1}$ that $\l\alpha_1,...,\alpha_n\r\in B$, hence $A_1\times....\times A_n\subseteq B$. $\blacksquare$

Every $X\in Ex(p)$ defines an ultrafilter $$\vec{U}(X,p)=\prod_{i=1}^{n+1}\prod_{j=1}^{|X_i|}U(\kappa_i(p),x_{i,j})$$ Note that $X(p)\in \vec{U}(X,p)$ by the definition of the product. Fix an extension type $X$ of $p$, every extension of $p$ of type $X$ corresponds to some element in the set $X(p)$ which is just a product of large sets. 

Let us state here some combinatorial properties, the proof can be found in \cite{TomMoti}. 
 
\begin{lemma}\label{stab}
 Let $\kappa_{1} \leq \kappa_{2} \leq ... \leq \kappa_{n}$ be a non descending  finite sequence of measurable cardinals and let $U_{1},...,U_{n}$ be normal measures\footnote{A measure over a measurable cardinal $\lambda$ is a $\lambda$-complete nonprincipal ultrafilter over $\lambda$.} over them respectively. Assume $
F:\overset{n}{\underset{i=1}{\prod }}A_{i}\longrightarrow \nu$ where $ \nu <\kappa_1$ and $A_{i} \in U_{i}$. Then there exists $H_{i}\subseteq A_{i}$,  $H_{i} \in U_{i}$ such that $ \overset{n}{\underset{i=1}{\prod}}H_{i}$ is homogeneous for $F$ i.e. $|Im(F\restriction \overset{n}{\underset{i=1}{\prod}}H_{i})|=1$.
\end{lemma}$\blacksquare$

Let $F:\prod_{i=1}^nA_i\rightarrow X$ be a function, and $I\subseteq\{1,...,n\}$. Let $$(\prod_{i=1}^nA_i)_I=\{\vec{\alpha}\restriction I\mid \vec{\alpha}\in\prod_{i=1}^nA_i\}$$
For $\vec{\alpha}'\in(\prod_{i=1}^nA_i)_I$, define $F_I(\vec{\alpha}')=F(\vec{\alpha})$  where $\vec{\alpha}\restriction I=\vec{\alpha}'$. With no further assumption, $F_I$ is not a well defined function.
 
\begin{lemma}\label{ImpIn}
 Let $\kappa_{1} \leq \kappa_{2} \leq ... \leq \kappa_{n}$ be a non descending  finite sequence of measurable cardinals and let $U_{1},...,U_{n}$ be normal measures over them respectively. Assume $
F:\overset{n}{\underset{i=1}{\prod }}A_{i}\longrightarrow B$ where $B$ is any set, and $A_{i} \in U_{i}$. Then there exists $H_{i}\subseteq A_{i}$, $H_{i} \in U_{i}$ and set $I \subseteq  \{1,...,n\}$  such that $F_I\restriction(\overset{n}{\underset{i=1}{\prod}}H_{i})_I:(\overset{n}{\underset{i=1}{\prod}}H_{i})_I\rightarrow B$ is well defined  and injective. 
\end{lemma}
\begin{definition}
Let $F:\prod_{i=1}^nA_i\rightarrow X$ be a function. An important coordinate is an index $r\in\{1,...,n\}$, such that for every $\vec{\alpha},\vec{\beta}\in \prod_{i=1}^nA_i$, $F(\vec{\alpha})=F(\vec{\beta})\rightarrow \vec{\alpha}(r)=\vec{\beta}(r)$.
\end{definition}

Lemma \ref{ImpIn} ensures the existence of a set $I$ of important coordinates, such that $I$ is ideal in the sense that removing any coordinate defect definition of $F_I$ as a function, and any coordinate  outside of $I$ is redundant.

We will need here another property that does not appear in \cite{TomMoti}.
\begin{lemma}\label{SepIm}
Let $\kappa_1\leq\kappa_2\leq...\leq\kappa_n$ and $\theta_1\leq\theta_2...\leq\theta_m$ be a non descending  finite sequences of measurable cardinals with corresponding normal measures
$U_1,....,U_n, W_1,...,W_m$.
Let $$F:\prod_{i=1}^nA_i\rightarrow X, \ G:\prod_{j=1}^mB_j\rightarrow X$$ be functions such that $X$ is any set, $A_i\in U_i$ and $B_j\in W_j$. Assume that $I\subseteq\{1,...,n\}$ and $ J\subseteq\{1,...,m\}$ are sets of important coordinates for $F,G$ respectively obtained by lemma \ref{ImpIn}. Then
there exists $A'_i\in U_i$ and $B'_j\in W_j$. such that one of the following holds \begin{enumerate}
    \item $Im(F\restriction\prod_{i=1}^nA'_i)\cap Im(G\restriction\prod_{j=1}^mB'_j)=\emptyset$.
    \item $(\prod_{i=1}^nA'_j)_I=(\prod_{j=1}^mB'_j)_J$ and $F_I\restriction(\prod_{i=1}^nA_i')_I=G_J\restriction(\prod_{j=1}^mB_j')_J$.
\end{enumerate}

\end{lemma}
\pr Fix $F,G$. Let us first deal with some trivial cases:
If $I=J=\emptyset$ i.e. $F,G$ are constantly $d_F,d_G$, respectively. Either $d_1\neq d_2$ and $(1)$ holds, or $d_1=d_2$ and $(2)$ holds. If $I=\emptyset$ and $j_0\in J\neq\emptyset$, then $F$ constantly $d_F$.
If $d_F\notin Im(G)$ then $(1)$ holds, otherwise, there is $\vec{\beta}$ such that $G(\vec{\beta})=d_F$,  remove $\vec{\beta}_{j_0}$
from $B_{j_0}$, then . If $\vec{\beta}'\in B_1\times...\times B_{j_0}\setminus\{\vec{\beta}_{j_0}\}\times...\times B_m$, then $G(\vec{\beta}')\neq d_F$, otherwise, $\vec{\beta}'\restriction J=\vec{\beta}\restriction J$
and in particular $\vec{\beta}_{j_0}=\vec{\beta}'_{j_0}$, contradiction. Similarly, if $J=\emptyset$ and $I\neq \emptyset$ then we can ensure $(1)$. We assume that $I,J\neq\emptyset$, also, without loss of generality, assume that $\kappa_1\leq\theta_1$. If $\kappa_1<\theta_1$ shrink the sets so that $\min(B_1)>\kappa_1$. By induction on $\langle n,m\rangle\in \mathbb{N}_+^2$ with respect to the lexicographical order.

\textbf{Case 1: Assume that $n=m=1$},  Assume that $I,J\neq \emptyset$. Define $$H_1:A_1\times B_1\rightarrow\{0,1\}, \ \ H_1(\alpha,\beta)=1\Leftrightarrow F(\alpha)=G(\beta)$$ 
By \ref{stab}, shrink $A_1,B_1$ to $A_1',B_1'$ so that $H_1$ is constant with colors $c_1$. If $c_1=1$ by fixing $\alpha$ we see that $G$ is constant on $B_1'$ with some value $\gamma$. It follows that $J=\emptyset$, contradiction. Assume that $c_1=0$, then for every $\alpha\in A_1,\beta\in B_1$ if $\alpha<\beta$ then $H_1(\alpha,\beta)=0$, which implies $F(\alpha)\neq G(\beta)$. This suffices for the case $\kappa_1<\theta_1$.
If $\kappa_1=\theta_1$, then it is possible that $\beta<\alpha$, so define
$$H_2:B_1\times A_1\rightarrow\{0,1\} \ \ H_2(\beta,\alpha)=1\Leftrightarrow F(\alpha)=G(\beta)$$
Again shrink the sets so that $H_2$ is constantly $c_2\in\{0,1\}$. In case $c_2=1$ we reach a similar contradiction to $c_1=1$. Assume that $c_2=0$, together with $c_1=0$, it follows that if $\beta\neq\alpha$ then $F(\alpha)\neq G(\beta)$. If $U_1\neq W_1$ then we can avoid the situation where $\alpha=\beta$ by separating $A_1',B_1'$  and
conclude that $$Im(F\restriction A_1')\cap Im(G\restriction B_1')=\emptyset$$ If $U_1=W_1$ then define $$H_3: A_1'\cap B_1'\rightarrow \{0,1\}, \ \ H_3(\alpha)=1\Leftrightarrow F(\alpha)=G(\alpha)$$
Again by \ref{stab} we can assume that $H_3$ is constant on $A^*$, if that constant is $1$ then we have $F\restriction A^*=G\restriction A^*$ (in particular $I=J=\{1\}$ and $F_I\restriction (A^*)_I=G_J\restriction (A^*)_J$) otherwise, $$Im(F\restriction A^*)\cap Im(G\restriction A^*)=\emptyset$$

\textbf{Case 2a: Assume $n=1$ and $m>1$}, By the assumption that $I,J\neq\emptyset$, $I=\{1\}$. Define $$H_1:A_1\times\prod_{j=1}^{m}B_j\rightarrow 
\{0,1\}, \ \ H_1(\alpha,\vec{\beta})=1\Leftrightarrow F(\alpha)=G(\vec{\beta})$$
Shrink the sets so that $H_1$ is constantly $c_1$. As before, if $c_1=1$ then $F,G$ are constant which is a contradiction.  Assume that $c_1=0$, which means that whenever $\alpha<\beta_1$, then $F(\alpha)\neq G(\vec{\beta})$. As before, if $\kappa_1<\theta_1$ then we are done. If $\kappa_1=\theta_1$,
for each $\beta\in B_1$, consider the function $$G_{\beta}:\prod_{j=2}^mB_j\setminus(\beta+1)\rightarrow X, \ G_\beta(\vec{\beta})=G(\beta^{\smallfrown}\vec{\beta})$$
Apply induction to $F$ and $G_\beta$, $\{1\},J\setminus\{1\}$ to find
$$A_1^{\beta}\in U_1, \ \ B^\beta_j\in W_j\text{ for }2\leq j\leq m$$ such that one of the following holds:
\begin{enumerate}
    \item  $A^{\beta}_1=(\prod_{j=1}^mB^{\beta}_j)_{J\setminus\{1\}}$, and $F\restriction A^{\beta}_1=(G_\beta)_{J\setminus\{1\}}\restriction(\prod_{j=2}^mB^{\beta}_j)_{J\setminus\{1\}}$.
    \item $Im(F\restriction A^{\beta}_1)\cap Im(G_{\beta}\restriction\prod_{j=2}^mB^{\beta}_j)=\emptyset$.
\end{enumerate}
Denote by $j_\beta\in\{1,2\}$ the relevant case. There is $B_1'\subseteq B_1$,  $B_1'\in W_1$, and $j^*\in\{1,2\}$ such that for every $\beta\in B_1'$, $j_\beta=j^*$. Let
$$A'_1=\underset{\beta\in B'_1}{\Delta}A^{\beta}_1, \ B'_j=\underset{\beta\in B'_1}{\Delta}B^{\beta}_j\text{
(Since }\theta_1=\kappa_1\text{ we can take the diagonal intersection)}$$
If $j^*=1$, then since $A^{\beta}_1=(\prod_{j=1}^mB^{\beta}_j)_{J\setminus\{1\}}$, it follows that $J=\{j_0\}$ and $A^{\beta}_1=B^{\beta}_{j_0}$ thus $A'_1=B'_{j_0}$. Also for $\beta_1,\beta'_1$, and some $\beta_1,\beta'_1<\beta_2,...,\beta_m$ in the product, $$G(\l\beta_1,...,\beta_m\r)=(mb nG_{\beta_1})_{j_0}(\beta_{j_0})=F(\beta_{j_0})=(G_{\beta_1'})_{j_0}(\beta_{j_0})=G(\l\beta_1',..,\beta_n\r)$$
Hence $1\notin J$, $A'_1=B'_{j_0}=(\prod_{j=1}^mB'_j)_J$ and $F_{1}\restriction A'_1=G_{j_0}\restriction B'_{j_0}$.

If $j^*=2$, for every $\l\beta_1,...\beta_m\r\in\prod_{j=1}^m B'_j$, $G(\l\beta_1,...,\beta_m\r)\in Im(G_{\beta_1}\restriction \prod_{j=1}^m B^{\beta}_j)$. Now if $\beta_1<\alpha\in A'_1$ then by definition of diagonal intersection $\alpha\in A^{\beta_1}_1$ and therefore, $F(\alpha)\in Im(F\restriction A^{\beta_1}_1)$ and we are done. Together with the assumption that $c_1=0$, we conclude that if $\alpha\neq\beta_1$ then $F(\alpha)\neq G(\vec{\beta})$. As before, we can avoid this situation if $U_1\neq W_1$. Assume that $U_1=W_1$, and assume that $A'_1=B'_1$.
 Let $$T_1:A'_1\times\prod_{j=2}^m B_j'\rightarrow\{0,1\}, \ \ T_1(\alpha,\vec{\beta})=1\Leftrightarrow
F(\alpha)=G(\alpha,\vec{\beta})$$ We shrink $A_1'$ and $B_j'$ so that $T_1$ is constantly $d_1$. If $d_1=0$ then we have eliminated the possibility of $\alpha=\beta$, and again we conclude that 
$$Im(F\restriction\prod_{i=1}^nA_i')\cap Im(G\restriction \prod_{j=1}^mB'_j)=\emptyset$$

If $d_1=1$ then $G$ only depends on $B'_1$ i.e. $J=\{1\}$, hence $(\prod_{j=1}^m B_j')_{\{1\}}=A'_1$ and $F\restriction A_1'=G_{\{1\}}\restriction A_1'$. 

\textbf{Case 2b: Assume $n>1$ and $m=1$}  Then by the assumption that $I,J\neq\emptyset$, it follows that $J=\{1\}$.
For $\alpha\in A_1$ define the functions $$F_\alpha:\prod_{i=2}^nA_i\setminus(\alpha+1)\rightarrow X, \ \ F_\alpha(\vec{\alpha})=F(\alpha,\vec{\alpha})$$ By the induction hypothesis applied to $F_\alpha,G$  and $I\setminus \{1\},\{1\}$, we obtain $$A^{\alpha}_i\in U_i \text{ for }  2\leq i\leq n, \ \ B^{\alpha}_j\in W_j \text{ for }1\leq j\leq m$$ such that one of the following holds:
\begin{enumerate}
    \item  $(\prod_{i=2}^nA^{\alpha}_i)_{I\setminus\{1\}}=B^{\alpha}_1$, and $(F_\alpha)_{I\setminus\{1\}}\restriction(\prod_{i=2}^nA^{\alpha}_i)_{I\setminus\{1\}}=G\restriction B^{\alpha}_1$.
    \item $Im(F_\alpha\restriction\prod_{i=2}^nA^{\alpha}_i)\cap Im(G\restriction B^{\alpha}_1)=\emptyset$.
\end{enumerate}
Denote by $i_\alpha\in\{1,2\}$ the relevant case. There is $A_1'\subseteq A_1$,  $A_1'\in U_1$, and $i^*\in\{1,2\}$ such that for every $\alpha\in A_1'$, $i_\alpha=i^*$. Let
$$A'_i=\underset{\alpha\in A_1}{\Delta}A^{\alpha}_i, \ B'_1=\underset{\alpha\in A_1}{\Delta}B^{\alpha}_1\text{
(Since }\theta_1\geq\kappa_1\text{ we can take the diagonal intersection)}$$ If $i^*=1$, then $(\prod_{i=2}^nA^{\alpha}_i)_{I\setminus\{1\}}=B^{\alpha}_1$, hence by $I=\{i_0\}$. Note that, $A^{\alpha}_{i_0}=B^{\beta}_{1}$ and in turn it follows that, $A'_{i_0}=B'_{1}\in U_{i_0}\cap W_{1}$.

Let $\alpha,\alpha'\in A'_1$, and $\alpha_1,\alpha'_1<\alpha_2<...<\alpha_n$ in the product, then $$F(\l\alpha_1...\alpha_n\r)=(F_{\alpha_1})_{\{i_0\}}(\alpha_{i_0})=G(\alpha_{i_0})=(F_{\alpha_1'})_{\{i_0\}}(\alpha_{i_0})=F(\l\alpha'_1...\alpha_n\r)$$ From this it follows that $1\notin I$, $B'_1=A'_{i_0}=(\prod_{i=1}^nA'_i)_I$ and $F_{I}\restriction A'_{i_0}=G\restriction B'_1$. Assume $i^*=2$ which means that for every $\l\alpha_1,...,\alpha_n\r\in \prod_{i=1}^nA'_1$, by definition of diagonal intersection, $\l\alpha_2,...,\alpha_n\r\in \prod_{i=2}^nA^{\alpha_1}_i$ hence
$$F(\l\alpha_1,...,\alpha_n\r)=F_{\alpha_1}(\l \alpha_2,..,\alpha_n\r)\in Im(F_{\alpha_1}\restriction \prod_{i=2}^nA^{\alpha_1}_i)$$
If $\beta\in  B'_1$, we cannot conclude automatically that $\beta\in B^{\alpha_1}_1$, since it is possible that $\beta_1\leq \alpha_1$. If $\kappa_1<\theta_1$, then $\beta_1\leq \alpha_1$ is impossible, thus, $\beta\in B^{\alpha_1}_1$ and
$G(\beta_1)\in Im(G\restriction B^{\alpha_1}_1)$.
since $i_{\alpha_1}=i^*=2$, it follows that $F(\l\alpha_1,...,\alpha_n\r)\neq G(\beta_1)$ which implies $$Im(F\restriction \prod_{i=1}^nA'_i)\cap Im(G\restriction B'_1)=\emptyset$$

 If $\theta_1=\kappa_1$, then we define 
 $$H_2:B_1\times\prod_{i=1}^{n}A_i\rightarrow 
\{0,1\}, \ \ H_2(\beta,\vec{\alpha})=1\Leftrightarrow F(\vec{\alpha})=G(\beta)$$
Shrink the sets so that $H_2$ is constantly $c_1$. As before, if $c_1=1$ then $F,G$ are constant which is a contradiction. Assume that $c_1=0$, which means that whenever $\beta<\alpha_1$, then $F(\vec{\alpha})\neq G(\beta)$. So we are left with the case $\alpha_1=\beta$, if $U_1\neq W_1$ then we can eliminate such example, and if $U_1=W_1$, consider $A^*_1=A'_1\cap B'_1$.
$$T_2:A^*_1\times\prod_{i=2}^n A_i'\rightarrow\{0,1\}, \ \ T_2(\alpha,\vec{\alpha})=1\Leftrightarrow
G(\alpha)=F(\alpha,\vec{\alpha})$$ We shrink $A_1^*$ and $A_i'$ so that $T_2$ is constantly $d_1$. If $d_1=0$ then we have eliminated the possibility of $\alpha=\beta$, and again we conclude that 
$$Im(F\restriction\prod_{i=1}^nA_i')\cap Im(G\restriction A^*_1)=\emptyset$$

If $d_1=1$ then $F$ only depends on $A^*_1$ i.e. $I=\{1\}$, hence $(A^*_1\times\prod_{i=2}^n A_i')_{\{1\}}=A^*_1$ and $G\restriction A^*_1=G_{\{1\}}\restriction A^*_1$. 

\textbf{Case 3: Assume $n,m>1$}
 
For $\alpha\in A_1$ define the functions $$F_\alpha:\prod_{i=2}^nA_i\setminus(\alpha+1)\rightarrow X, \ \ F_\alpha(\vec{\alpha})=F(\alpha,\vec{\alpha})$$ By the induction hypothesis applied to $F_\alpha,G$  and $I\setminus \{1\},J$, we obtain $$A^{\alpha}_i\in U_i \text{ for }  2\leq i\leq n, \ \ B^{\alpha}_j\in W_j \text{ for }1\leq j\leq m$$ such that one of the following holds:
\begin{enumerate}
    \item  $(\prod_{i=2}^nA^{\alpha}_i)_{I\setminus\{1\}}=(\prod_{j=1}^mB^{\alpha}_j)_{J}$, and $(F_\alpha)_{I\setminus\{1\}}\restriction(\prod_{i=2}^nA^{\alpha}_i)_{I\setminus\{1\}}=G_J\restriction(\prod_{j=1}^mB^{\alpha}_j)_{J}$.
    \item $Im(F_\alpha\restriction\prod_{i=2}^nA^{\alpha}_i)\cap Im(G\restriction\prod_{j=1}^mB^{\alpha}_j)=\emptyset$.
\end{enumerate}
Denote by $i_\alpha\in\{1,2\}$ the relevant case. There is $A_1'\subseteq A_1$,  $A_1'\in U_1$, and $i^*\in\{1,2\}$ such that for every $\alpha\in A_1'$, $i_\alpha=i^*$. Let
$$A'_i=\underset{\alpha\in A_1}{\Delta}A^{\alpha}_i, \ B'_j=\underset{\alpha\in A_1}{\Delta}B^{\alpha}_j\text{
(Since }\theta_1\geq\kappa_1\text{ we can take the diagonal intersection)}$$ If $i^*=1$, then $(\prod_{i=2}^nA^{\alpha}_i)_{I\setminus\{1\}}=(\prod_{j=1}^mB^{\alpha}_j)_{J}$, denote by $I\setminus\{1\}=\{i_1,...,i_k\}$, $J=\{j_1,...,j_k\}$. Note that for every $1\leq r\leq k$, $A^{\alpha}_{i_r}=B^{\beta}_{j_r}$, thus $A'_{i_r}=B'_{j_r}\in U_{i_r}\cap W_{j_r}$.  It follows that,
$$(\prod_{i=1}^nA'_i)_{ I\setminus\{1\}}=(\prod_{j=1}^mB'_j)_{J}$$
Let $\alpha,\alpha'\in A'_1$, $\vec{\alpha}\in \prod_{i=2}^nA_i'$ with $\min(\vec{\alpha})>\alpha,\alpha'$, then  $$F_\alpha(\vec{\alpha})=(F_\alpha)_{I\setminus\{1\}}(\vec{\alpha}\restriction I)=G_{J}(\vec{\alpha}\restriction I)=(F_{\alpha'})_{I\setminus\{1\}}(\vec{\alpha}\restriction I)=F_{\alpha'}(\vec{\alpha})$$ From this it follows that $1\notin I$ and $F_{I}=F_{I\setminus\{1\}}=G_{J}$. Assume $i^*=2$ which means that for every $\l\alpha_1,...,\alpha_n\r\in \prod_{i=1}^nA'_1$, by definition of diagonal intersection, $\l\alpha_2,...,\alpha_n\r\in \prod_{i=2}^nA^{\alpha_1}_i$ hence
$$F(\l\alpha_1,...,\alpha_n\r)=F_{\alpha_1}(\l \alpha_2,..,\alpha_n\r)\in Im(F_{\alpha_1}\restriction \prod_{i=2}^nA^{\alpha_1}_i)$$
If $\vec{\beta}\in\prod_{j=1}^mB'_j$, we cannot conclude automatically that $\vec{\beta}\in \prod_{j=1}^mB^{\alpha_1}_j$, since it is possible that $\beta_1\leq \alpha_1$. If $\kappa_1<\theta_1$, then $\beta_1\leq \alpha_1$ is impossible, thus, $\vec{\beta}\in \prod_{j=1}^mB^{\alpha_1}_j$ and
$G(\l\beta_1,...,\beta_n\r)\in Im(G\restriction \prod_{j=1}^nB^{\alpha_1}_j)$.
since $i_{\alpha_1}=i^*=2$, it follows that $F(\l\alpha_1,...,\alpha_n\r)\neq G(\l\beta_1,...,\beta_n\r)$ which implies $$Im(F\restriction \prod_{i=1}^nA'_i)\cap Im(G\restriction \prod_{j=1}^nB'_j)=\emptyset$$

 If $\theta_1=\kappa_1$, we repeat the same process, we use $G_\beta$ and fix $F$, denoting $j_{\beta}$ the relevant case, shrink the sets so that $j^*$ is constant. In case $j^*=1$ the proof is the same as $i^*=1$. So we assume that $i^*=j^*=2$, meaning that for every $\langle\alpha\r^{\smallfrown}\vec{\alpha}\in \prod_{i=1}^nA_i'$ and every $\l\beta\r^{\smallfrown}\vec{\beta}\in\prod_{j=1}^mB'_j$
$$\alpha\neq\beta\rightarrow F(\alpha,\vec{\alpha})\neq G(\beta,\vec{\beta})$$

We are left with the case $\alpha=\beta$.

\textbf{Case 3a: Assume that $U_1\neq W_1$}   Then we can just shrink the sets $A'_1,B'_1$ so that $A_1'\cap B_1'=\emptyset$. Together with the construction of case $3$, conclude that
$$Im(F\restriction\prod_{i=1}^nA_i')\cap Im(G\restriction \prod_{j=1}^mB'_j)=\emptyset$$

\textbf{Case 3b: Assume that $U_1=W_1$}, then we shrink the sets so that $A'_1=B'_1$.  For every $\alpha\in A'_1$ we apply the induction hypothesis to the functions $F_\alpha,G_\alpha$, this time denoting the cases by $r^*$. If $r^*=2$, then we have eliminated the possibility of $F(\alpha,\vec{\alpha})=G(\alpha,\vec{\beta})$, together with $i^*=2,j^*=2$ we are done. Finally, assume $r^*=1$, namely that for $$I^*:=I\setminus\{1\}\subseteq\{2,...,n\},\  J^*:=J\setminus\{1\}\subseteq\{2,...,m\}$$
 We have
 $$ (\prod_{i=2}^nA'_i)_{I^*}=(\prod_{j=2}^mB'_j)_{J^*}\text{ and }(F_\alpha)_{I^*}\restriction(\prod_{i=2}^nA'_i)_{I^*}=(G_{\alpha})_{J^*}\restriction(\prod_{j=2}^mB'_j)_{J^*}$$  Since $A'_1=B'_1$ it follows that $$(*) \ \ \ \ \ (\prod_{i=1}^nA'_i)_{ I^*\cup\{1\}}=(\prod_{j=1}^mB'_j)_{\in J^*\cup\{1\}}\text{ and }(F_\alpha)_{I^*\cup\{1\}}\restriction(\prod_{i=2}^nA'_i)_{I^*}=(G_{\alpha})_{J^*}\restriction(\prod_{j=2}^mB'_j)_{J^*\cup\{1\}}$$ Since if $\l\alpha\r^{\smallfrown}\vec{\alpha}\in (\prod_{i=1}^nA'_i)_I$, $$F_{I^*\cup\{1\}}(\alpha,\vec{\alpha})=(F_\alpha)_{I^*}(\vec{\alpha})=(G_\alpha)_{J^*}(\vec{\alpha})=G_{J^*\cup\{1\}}(\alpha,\vec{\alpha})$$
We claim that $1\in I$ if and only if $1\in J$. By symmetry, it suffices to prove one implication,  for example, if $1\in I$, then $I=I^*\cup\{1\}$, take $\vec{\alpha}\restriction I,\vec{\alpha}'\restriction I\in (\prod_{i=1}^nA'_i)_I$ which differs only at the first coordinate, therefore $F(\vec{\alpha})\neq F(\vec{\alpha}')$. By $(*)$, there are $\vec{\beta},\vec{\beta}'\in \prod_{i=1}^m B'_i$ such that $$\vec{\beta}\restriction (J^*\cup\{1\})=\vec{\alpha}\restriction I\text{ and }\vec{\beta}'\restriction (J^*\cup\{1\})=\vec{\alpha}'\restriction I$$ It follows that from $(*)$ that $G(\vec{\beta})=F(\vec{\alpha})\neq F(\vec{\alpha}')=G(\vec{\beta}')$, therefore $1\in J$. 

In any case, $F_I\restriction(\prod_{i=1}^nA'_i)_I =G_J\restriction(\prod_{i=1}^mB'_i)_J$. 
$\blacksquare$
\section{ The main result}
Let us turn to prove the main result (theorem \ref{MainResaultPartone}) for Magidor forcing with $o^{\vec{U}}(\kappa)<\kappa$.
The proof presented here is based on what was done in \cite{TomMoti} and before that in \cite{PrikryCaseGitikKanKoe}, it is a proof by induction of $\kappa$. 
\subsection{Short Sequences}
In this section we prove the theorem for sets $A$ of small cardinality. 
\begin{proposition}\label{amalgamate1} Let $p\in\Mfor$ be any condition, $X$ an extension type of $p$. For every $\vec{\alpha}\in X(p)$ let $p_{\vec{\alpha}}\geq^*p^{\frown}\vec{\alpha}$. Then there exists $p\leq^*p^*$ such that for every $\vec{\beta}\in X(p^*)$, every $p^{*\frown}\vec{\beta}\leq q$ is compatible with $p_{\vec{\beta}}$ .
\end{proposition}
\pr By induction of $l(X)$. If $l(X)=1$, $X=\langle \xi\rangle$, then $\vec{U}(X,p)=U(\kappa_i(p),\xi)$ and $X(p)=B_{i,\xi}(p)$. For each $\beta\in B_{i,\xi}(p)$
$$p_\beta=\langle\langle\kappa_1(p),A^{\beta}_1\rangle,...,\langle\kappa_{i-1}(p),A^{\beta}_{i-1}\rangle,\langle \beta,B_\beta\rangle,\langle\kappa_i(p),A^{\beta}_i\rangle,...,\langle\kappa,A_{\beta}\rangle\rangle$$
For $j>i$ let $A^*_j=\cap_{\beta\in B_{i,\xi}(p)}A^{\beta}_j$.
For $j<i$ we can find $A^*_j$ and shrink $B_{i,\xi}(p)$ to $E_\xi$ so that for every $\beta\in E_{\xi}$ and $j<i$ $A^{\beta}_j=A^*_j$.
For $i$, first let $E=\Delta_{\alpha\in B_{i,\xi}(p)} A^{\beta}_i$.
By ineffability of $\kappa_i(p)$ we can find $A^*_\xi\subseteq E_\xi$ and a set $B^*\subseteq\kappa_i(p)$ such that for every $\beta\in A^*_\xi$ $B^*\cap\beta=B_\beta$.
Claim that $B^*\in U(\kappa_i(p),\gamma)$ for every $\gamma<\xi$, $$Ult(V,U(\kappa_i(p),\xi))\models B^*=j_{U(\kappa_i(p),j)}(B^*)\cap\kappa_i(p)$$ and since
$$\{\beta<\kappa\mid B^*\cap\beta\in\cap\vec{U}(\beta)\}\in U(\kappa_i(p),\xi)$$
it follows that $B^*\in\cap j_{U(\kappa_i(p),\xi)}(\vec{U})(\kappa_i(p))$.
By coherency $B^*\in\cap_{\gamma<\xi}U(\kappa_i(p),\gamma)$.
Define  $$A^*_i=B^*\uplus A^*_{\xi}\uplus(\underset{\xi<i}{\cup}E_i)\in \cap\vec{U}(\kappa_i(p))$$
Let $q\geq p^{*\frown}\beta$ and suppose that $q\geq^*(p^{*\frown}\beta)^{\frown}\vec{\gamma}$. Then every $\gamma\in\vec{\gamma}$ such that $\gamma>\beta$ belong to some $A^*_j\setminus\beta$ for $j\geq i$, and by the definition of these sets $\gamma\in A^{\beta}_j$. If $\gamma<\kappa_{i-1}$ then also $\gamma\in A^*_j$ for some $j<i$.  Since $\beta\in E_\xi$ it follows that $A^{\beta}_j=A^*_j$ so $\gamma\in A^{\beta}_j$. For $\gamma\in(\kappa_{i-1},\beta)$, by definition of the order we have $o^{\vec{U}}(\gamma)<o^{\vec{U}}(\beta)=\xi$ and therefore $\gamma\in A^*_{i,\eta}\cap\beta$ for some $\eta<\xi$, but $$A^*_{i,\eta}\cap\beta\subseteq B^*\cap\beta=B_\beta$$ it follows that $q,p_\beta$ are compatible. For general $X$, fix $\min(\vec{\beta})=\beta$. Apply the induction hypothesis to $p^{\frown}\beta$ and $p_{\vec{\beta}}$ to find $p^*_\beta\geq^* p^{\frown}\beta$. Next apply the case $n=1$ to $p^*_\beta$ and $p$, find $p^*\geq p$. Let $q\geq p^{*\frown}\vec{\beta}$ and denote $\beta=\min(\vec{\beta})$ then $q$ is compatible with  $p^*_\beta$ thus let $q'\geq q, p^*_\beta$. Since $q'\geq p^*_\beta$ and $q'\geq p^{*\frown}\vec{\beta}$ it follows that $q'\geq p^{*\frown}_\beta\vec{\beta}$. Therefore there is $q''\geq q', p_{\vec{\beta}}$.
$\blacksquare$
\begin{lemma}\label{exT}
Let $\lambda<\kappa$, $p\in\Mfor\restriction(\lambda,\kappa)$, $q\in \Mfor\restriction\lambda$ and $X\in Ex(p)$. Also. let $\underaccent{\sim}{x}$ be an ordinal $\Mfor$-name. There is $p\leq^*p^*$ such that
    $$\text{ If } \exists \vec{\alpha}\in X(p^*) \ \exists p'\geq^*p^{*\frown}\vec{\alpha}\ \langle q,p'\rangle|| \ \underaccent{\sim}{x}\ \text{ Then } \forall \vec{\alpha}\in X(p^*) \langle q,p^{*\frown}\vec{\alpha}\rangle|| \underaccent{\sim}{x}$$
\end{lemma}
\pr Fix $p,\lambda,q,X$ as in the lemma. Consider the set
$$B_0=\{\vec{\beta}\in X(p)\mid \exists p'\prescript{*}{}{\geq}p^{\frown}\vec{\beta} \ s.t. \ \langle q,p'\rangle||\underaccent{\sim}{x  }\}$$
One and only one of $B_0$ and $X(p)\setminus B_0$ is in $\vec{U}(X,P)$. Denote this set by $A'$. By proposition \ref{Normgenerate}, we can find $A'_{i,j}\in U(\alpha_i,x_{i,j})$ such that $\prod_{i=1}^{l(p)+1}\prod_{j=1}^{|X_i|}A'_{i,j}\subseteq A'$, let $p\leq^*p'$ be the condition obtained by shrinking $B_{i,j}(p)$ to $A'_{i,j}$ so that $X(p')=\prod_{i=1}^{n+1}\prod_{j=1}^{|X_i|}A'_{i,j}$. If $$\exists \vec{\beta}\in X(p') \ \exists p''\prescript{*}{}{\geq} p'^{\frown}\vec{\beta}\ \langle q,p''\rangle|| \ \underaccent{\sim}{x}$$
Then $\vec{\beta}\in B_0\cap A'$ and therefore $B_0=A'$ , we conclude that 
$$\forall \vec{\beta}\in X(p') \ \exists p_{\vec{\beta}}\prescript{*}{}{\geq} p'^{\frown}\vec{\beta}\ \langle q,p_{\vec{\beta}}\rangle|| \ \underaccent{\sim}{x}$$
By proposition \ref{amalgamate1} we can amalgamate all these $p_{\vec{\beta}}$ to find $p'\leq^* p^*$ such that for every $\vec{\beta}\in X(p^*)$, $p^{*\frown}\vec{\beta}$ decides $\underaccent{\sim}{x}$, then $p^*$ is as wanted.
$\blacksquare$
\begin{lemma}\label{Rad-prop} Consider the decomposition of \ref{dec} at some $\lambda\geq o^{\vec{U}}(\kappa)$ and let $\underaccent{\sim}{x}$ be a $\Mfor$-name for an ordinal. Then for every $p\in \Mfor\restriction(\lambda,\kappa)$, there exists $p\leq^*p^*$ such that for every $X\in Ex(p)$ and $q\in\Mfor\restriction\lambda$ the following holds:
$$If \ \exists \vec{\alpha}\in X(p^*) \ \exists p'\geq^*p^{*\frown}\vec{\alpha}\ \langle q,p'\rangle|| \ \underaccent{\sim}{x}\ Then \ \forall \vec{\alpha}\in X(p^*) \  \langle q,p^{*\frown}\vec{\alpha}\rangle|| \underaccent{\sim}{x}$$
\end{lemma}
\pr Fix $q\in\Mfor\restriction \lambda$ and and $X\in Ex(p)$. Use \ref{exT}, to find $p\leq^*p_{q,X}$ such that 
   $$If \ \exists \vec{\alpha}\in X(p_{q,X}) \ \exists p'\geq^*(p_{q,X})^{\frown}\vec{\alpha} \ s.t. \  \langle q,p'\rangle|| \ \underaccent{\sim}{x}\ Then \ \forall \vec{\alpha}\in X(p_{q,X}) \  \langle q,(p_{q,X})^{\frown}\vec{\alpha}\rangle|| \underaccent{\sim}{x}$$
By the definition of $\lambda$, the forcing $\Mfor\restriction(\lambda,\kappa)$ is $\leq^*$-$\max(|Ex(p)|^+,|\Mfor\restriction\lambda|^+)$-directed. Hence we can find $p\leq^* p^*$ so that for every $X,q$, $p_{q,X}\leq^* p^*$.
$\blacksquare$
\begin{lemma}\label{short} Let $A\in V[G]$ be a set of ordinals such that $|A|<\kappa$.
 Then there exists $C'\subseteq C_G$ such that $V[A]=V[C']$.

\end{lemma}
\pr Assume that $|A|=\lambda'<\kappa$ and let $\delta=max(\lambda',\otp(C_G))<\kappa$. Split $\Mfor$ as in proposition \ref{dec}. Find $p\in G$ such that some $\lambda\geq\delta$ appears in $p$. The generic $G$ also splits to $G=G_1\times G_2 $ where $G_1$ is the generic for Magidor forcing below $\lambda$ and by remark \ref{lift}, $G_2$ is $V[G_1]$-generic for the upper part of the forcing. Let   $\langle \underaccent{\sim}{a}_i \mid i<\lambda' \rangle$ be a $\Mfor$-name for $A$ in $V$ and $p\in\Mfor\restriction(\lambda,\kappa)$.  For every $i<\lambda'$ find $p\leq^* p_{i}$ as in lemma \ref{Rad-prop}, such that for every $q\in \Mfor\restriction \lambda$ and $X\in Ex(p)$ we have:
$$If \ \exists \vec{\alpha}\in X(p_i)\ \exists p_{i}^{\frown}\vec{\alpha}\leq^*p'  \  \langle q,p'\rangle \ || \ \underaccent{\sim}{a}_i\ Then \ \forall \vec{\alpha}\in X(p_i)\ \langle q,p_{i}^{\frown}\vec{\alpha}\rangle \ || \ \underaccent{\sim}{a}_i \ \ (*)$$
Since in $\Mfor\restriction(\lambda,\kappa)$ we have $\lambda^+$-closure for $\leq^*$, we can find a single $p_{i}\leq^* p_*$. Next, for every $i<\lambda'$, fix a maximal antichain $Z_i\subseteq \Mfor\restriction \lambda$ such that for every $q\in Z_i$ there is an extension type $X_{q,i}$ for which $\forall \vec{\alpha}\in p_*^{\frown}X_{q,i}\ \langle q,p_*^{\frown}\vec{\alpha}\rangle \ || \ \underaccent{\sim}{a}_i$, these antichains can be found using (*) and Zorn's lemma. Recall the sets $X_{q,i}(p_*)$ is a product of large sets. Define $F_{q,i}:X_{q,i}(p_*)\rightarrow On$ by
$$F_{q,i}(\vec{\alpha})=\gamma \ \ \ \Leftrightarrow \ \ \  \langle q,p_*^{\frown}\vec{\alpha}\rangle \Vdash \underaccent{\sim}{a}_i=\check{\gamma}$$
By lemma \ref{ImpIn}  we can assume that there are important coordinates $$I_{q,i}\subseteq\{1,...,Dom(X_{q,i}(p_{*}))\}$$ Fix $i<\lambda'$, for every $q,q'\in Z_i$ we apply lemma \ref{SepIm} to the functions $F_{q,i},F_{q,i'}$ and find $p_*\leq^*p_{q,q'}$ for which one of the following holds:
\begin{enumerate}
\item $Im(F_{q,i}\restriction A(X_{q,i},p_{q,q'}))\cap Im(F_{q',i}\restriction  A(X_{q',i},p_{q,q'}))=\emptyset$
\item 
 $(F_{q,i})_{I_{q,i}}\restriction( A(X_{q,i},p_{q,q'}))_{I_{q,i}}=(F_{q',i})_{I_{q',i}}\restriction (A(X_{q',i},p_{q,q'}))_{I_{q',i}}$
\end{enumerate}
Finally find $p^{*}$ such that for every $q,q'$, $p_{q,q'}\leq^*p^{*}$.
By density, there is such $p^{*}\in G_2$. We use $F_{q,i}$ to translate information from $C_G$ to $A$ and vice versa, distinguishing from \cite{TomMoti} this translation is made in $V[G_1]$ rather then $V$: For every $i<\lambda'$, $G_1\cap Z_i=\{q_i\}$. Use lemma \ref{Max}, to find $D_i\in X_{q_i,i}(p^{*})$ be such that $p^{*\frown}D_i\in G_2$, define $C_i=D_i\restriction I_{q_i,i}$ and let $C'=\underset{i<\lambda'}{\bigcup}C_i$.  Define as in \ref{Index}, $I(C_i,C')\in[\otp(C_G)]^{<\omega}$, since $\otp(C')\leq \otp(C_G)\leq\lambda$ and by proposition \ref{genericproperties}(6), $G_2$ does not add $\lambda$-sequences of ordinals below $\lambda$ to $V[G_1]$, we conclude that $\langle I(C_i,C')\mid i<\lambda'\rangle\in V[G_1]$. It follows that
$$(V[G_1])[A]=(V[G_1])[\langle C_i\mid i<\lambda'\rangle]=(V[G_1])[C']$$
In fact let us prove that $\langle C_i\mid i<\lambda'\rangle\in V[A]$. Indeed, define in $V[A]$ the sets $$M_i=\{q\in Z_i\mid a_i\in Im(F_{q,i})\}$$  then, for any $q,q'\in M_i$ $a_i\in Im(F_{q_i})\cap Im(F_{q',i})\neq\emptyset$. Hence $2$ must hold for $F_{q,i},F_{q',i}$ i.e.
 $$(F_{q,i})_{I_{q,i}}\restriction (X_{q,i}(p^{*}))_{I_{q,i}}=(F_{q',i})_{I_{q',i}}\restriction (X_{q',i}(p^{*}))_{I_{q',i}}$$
This means that no matter how we pick $q_i'\in M_i$, we will end up with the same function $(F_{q'_i,i})_{I_{q_i',i}}\restriction (X_{q'_i,i}(p^{*}))_{I_{q_i',i}}$. In $V[A]$, choose any $q'_i\in M_i$ and let $D_i'\in F_{q_i',i}^{-1}(a_i)$, $C_i'=D_i\restriction I_{q_i',i}$. Since $q_i,q'_i\in M_i$ we have $ C_i=C_i'$, hence $\langle C_i\mid i<\lambda'\rangle\in V[A]$. We still have to determine what information $A$ uses in the part of $G_1$, namely, $\{q'_i\mid i<\lambda'\},\langle I(C_i,C')\mid i<\lambda'\rangle\in V[A]$. This sets can be coded as a subset of ordinals below $(2^{\lambda})^+$, therefore, $$\{q'_i\mid i<\lambda'\},\langle I(C_i,C')\mid i<\lambda'\rangle\in V[G_1]$$ By the induction hypothesis applied to $G_1$, we can find $C''\subseteq C_{G_1}$ such that $$V[\{q'_i\mid i<\lambda'\},\langle I(C_i,C')\mid i<\lambda'\rangle]=V[C'']$$  Since all the information needed to restore $A$ is coded in $C'\uplus C''$, it is clear that $V[A]=V[C''\uplus C']$.
$\blacksquare$
\subsection{General Subsets of $\kappa$}
Assume that $A\in V[G]$ such that $A \subseteq \kappa$. For some $A$'s, the proof is similar to the one in \cite{TomMoti} works. This proof relies on the following lemma:
\begin{lemma}\label{suff}
Assume that $o^{\vec{U}}(\kappa)<\kappa$ and let $A \in V[G] ,\ \sup(A)=\kappa $.  Assume that $\exists C^*\subseteq C_G$ such that
\begin{enumerate}
 \item$C^*\in V[A]$ and
 $\forall\alpha<\kappa \  A\cap\alpha\in V[C^*]$
 \item $cf^{V[A]}(\kappa)<\kappa$
 \end{enumerate}
Then $ \exists C' \subseteq C_G $ such that $ V[A]=V[C']$.
\end{lemma}
\pr Let $\langle\alpha_i\mid i<\lambda\rangle\in V[A]$ be cofinal in $\kappa$. Since $|C^*|<\kappa$, by \ref{short}, we can find $C''\subseteq C_G$ such that $$V[C'']=V[C^*,\langle\alpha_i\mid i<\lambda\rangle]\subseteq V[A]$$ In $V[C'']$    choose for every $i$, a bijection $\pi_i:2^{\alpha_i}\rightarrow P^{V[C'']}(\alpha_i)$. Since $A\cap\alpha_i\in V[C'']$ there is $\delta_i$ such that $\pi_i(\delta_i)=A\cap\alpha_i$. Finally let $C'\subseteq C_G$ such that $$V[C']=V[C'',\langle \delta_i\mid i<\lambda\rangle]$$ We claim that $V[A]=V[C']$. Obviously, $C'\in V[A]$, for the other direction, $$\langle A\cap\alpha_i\mid i<\lambda\rangle=\langle\pi_i(\delta_i)\mid i<\lambda\rangle\in V[C']$$ Thus $A\in V[C']$.
$\blacksquare$
\begin{definition}
We say that $A\cap\alpha$ stabilizes, if $$\exists\alpha^*<\kappa. \ \forall\alpha<\kappa. \  A\cap\alpha\in V[A\cap\alpha^*]$$
\end{definition}

First we deal with $A$'s such that $A\cap\alpha$ does not stabilize.
\begin{lemma}
Assume $o^{\vec{U}}(\kappa)<\kappa$, $A\subseteq\kappa$ unbounded in $\kappa$ such that $A\cap\alpha$ does not stabilizes, then there is $C'\subseteq C_G$ such that  $V[C']=V[A]$.
\end{lemma}
\pr Work in $V[A]$, define the sequence $\langle\alpha_\xi \mid \xi<\theta\rangle$:
$$ \alpha_0=\min\{\alpha \mid V[A\cap\alpha]\supsetneq V\}$$
Assume that $\langle \alpha_\xi \mid \xi<\lambda\rangle$ has been defined and for every $\xi ,\ \alpha_\xi<\kappa$. 
If $\lambda=\xi+1$ then set
$$\alpha_\lambda=\min\{\alpha \mid V[A\cap\alpha]\supsetneq V[A\cap\alpha_\xi]\}$$
To see that $\alpha_\lambda$ is a well define ordinal below $\kappa$, note that by the assumption that $A$ does not stabilize, there is $\alpha<\kappa$ such that $A\cap\alpha\notin V[A\cap\alpha_{\xi}]$, hence $V[A\cap \alpha_\xi]\subsetneq V[A\cap \alpha]$.

If $\lambda$ is limit, define 
$$\alpha_\lambda=\sup(\alpha_\xi \mid \xi<\lambda)$$
 if $\alpha_\lambda=\kappa$ define $\theta=\lambda$ and stop.  The sequence $\langle \alpha_\xi \mid \xi<\theta\rangle\in V[A]$ is a continues, increasing unbounded sequence in $\kappa$. Therefore, $ cf^{V[A]}(\kappa)=cf^{V[A]}(\theta)$. Let us argue that $\theta<\kappa$. Work in $V[G]$, for every $\xi<\theta$ pick $C_\xi\subseteq C_G$ such that  $V[A\cap\alpha_\xi]=V[C_\xi]$. The map $\xi\mapsto C_\xi$ is injective from $\theta$ to  $P(C_G)$, by the definition of $\alpha_\xi$'s. Since $o^{\vec{U}}(\kappa)<\kappa$, $|C_G|<\kappa$, and $\kappa$ stays strong limit in the generic extension. Therefore
 $$\theta\leq|P(C_G)|=2^{|C_G|}<\kappa$$
Hence $\kappa$ changes cofinality in $V[A]$, according to lemma \ref{suff}, it remains to find $C^*$. Denote $\lambda=|C_G|$ and work in $V[A]$, for every $ \xi<\theta$,  $C_\xi\in V[A]$ (Although the sequence $\langle C_\xi \mid \xi<\theta\rangle$ may not be in $V[A]$). $C_\xi$ witnesses that
 $$ \exists d_\xi\subseteq \kappa. \ |d_{\xi}|\leq\lambda \text{  and } V[A\cap\alpha_\xi]=V[d_{\xi}]$$
Fix $d=\langle d_{\xi} | \xi<\theta\rangle\in V[A]$. It follows that $d$ can be coded as a subset of $\kappa$ of cardinality $\leq\lambda\cdot\theta<\kappa$. Finally, by \ref{short}, there exists $C^*\subseteq C_G$ such that $V[C^*]=V[d]\subseteq V[A]$ so $$\forall\alpha<\kappa. \ A\cap\alpha\in V[d_\xi]\subseteq V[C^*]$$
$\blacksquare$
\vskip 0.2 cm
Next we assume that $A\cap\alpha$ stabilizes on some $\alpha^*<\kappa$. By lemma \ref{short} There exists $C^*\subseteq C_G$ such that $V[A\cap\alpha^*]=V[C^*]$, if $A\in V[C^*]$ then we are done, assume that $A\notin V[C^*]$. To apply \ref{suff}, it remains to prove that $cf^{V[A]}(\kappa)<\kappa$. The subsequence $C^*$ must be bounded, denote $\kappa_1=\sup{(C^*)}<\kappa$ and $\kappa^*=\max(\kappa_1,\otp(C_G))$. Find $p\in G$ that decides the value of $\kappa^*$ and assume that $\kappa^*$ appear in $p$ (otherwise take some ordinal above it). As in lemma \ref{dec} we split
$$\Mfor/p\simeq\Big(\Mfor\restriction\kappa^*\Big)/\Big(p\restriction\kappa^*\Big)\times\Big(\Mfor\restriction(\kappa^*,\kappa)\Big)/\Big(p\restriction(\kappa^*,\kappa)\Big)$$
There is a complete subalgebra $\mathbb{P}$ of $RO(\Big(\Mfor\restriction\kappa^*\Big)/\Big(p\restriction\kappa^*\Big))$ such that $V[C^*]=V[H]$ for some $V$-generic filter  $H\subseteq\mathbb{P}$. Let $$\mathbb{Q}=\Big[\Big(\Mfor\restriction\kappa^*\Big)/\Big(p\restriction\kappa^*\Big)\Big]/C^*$$ be the quotient forcing completing $\mathbb{P}$ to
$\Big(\Mfor\restriction\kappa^*\Big)/\Big(p\restriction\kappa^*\Big)$. Finally note that $G$ is generic over $V[C^*]$ for $$\mathbb{S}=\mathbb{Q}\times\Big(\Mfor\restriction(\kappa^*,\kappa)\Big)/\Big(p\restriction(\kappa^*,\kappa)\Big)$$
\begin{lemma}
$cf^{V[A]}(\kappa)<\kappa$
\end{lemma}
\pr Let $G=G_1\times G_2$ be the decomposition such that $G_1$ is generic for $\mathbb{Q}$ above $V[C^*]$ and $G_2$ is $\Mfor\restriction(\kappa^*,\kappa)$ generic over $V[C^*][G_1]$. Let $\underaccent{\sim}{A}$ be a $\mathbb{S}$-name for $A$ in $V[C^*]$. and $\langle q_0,p_0\rangle\in G$ such that $$\langle q_0,p_0\rangle\Vdash \ `` \forall\alpha<\kappa \ \underaccent{\sim}{A}\cap\alpha \ is \ old" \  (i.e. \ in \ V[C^*])$$ Proceed by a density argument in $\Mfor\restriction(\kappa^*,\kappa))/p\restriction(\kappa^*,\kappa)$, let $p_0\leq p$, as in \ref{short} find $p\leq^*p^*$ such that for all $q_0\leq q\in\mathbb{Q}$ and $X\in Ex(p^*)$:
$$\exists\vec{\alpha}^{\smallfrown}\l\alpha\rangle\in X(p^*)\exists p'\geq^*p^{*\frown}\vec{\alpha}^{\smallfrown}\l\alpha\rangle \ \langle q,p'\rangle \ || \ \underaccent{\sim}{A}\cap\alpha\Rightarrow \forall \vec{\alpha}^{\smallfrown}\l\alpha\rangle\in X(p^*) . \ \langle q,p^{*\frown}\vec{\alpha}^{\smallfrown}\l\alpha \rangle\r \ || \ \underaccent{\sim}{A}\cap\alpha$$
Denote the consequent by $(*)_{X,q}$, since $\underaccent{\sim}{A}\cap\alpha$ is forced to be old, we will find many $q,X$ for which $(*)_{q,X}$ holds. For such $q,X$, for every $\vec{\alpha}^{\smallfrown}\l\alpha\r\in X(p^*)$ define the value forced for $\underaccent{\sim}{A}\cap\alpha$ by $a(q,\vec{\alpha},\alpha)$. Fix $q,X$ such that $(*)_{q,X}$ holds. Assume that the maximal measure which appears in $X$ is $U(\kappa_i(p),mc(X))$ and fix $\vec{\alpha}\in (X\setminus \{mc(X)\})(p^*)$. For every $\alpha\in B_{i,mc(X)}(p)\setminus\max(\vec{\alpha})$ the set
$a(q,\vec{\alpha},\alpha)\subseteq\alpha$ is defined. By ineffability, we can shrink $B_{i,mc(X)}(p)$ to $A^{q,\vec{\alpha}}_{i,mc(X)}$ and find a set $A(q,\vec{\alpha})\subseteq \kappa_i(p)$ such that for every $\alpha\in A^{q,\vec{\alpha}}_{i,mc(X)}$,
$A(q,\vec{\alpha})\cap\alpha=a(q,\vec{\alpha},\alpha)$ define $$A'_{i,mc(X)}=\underset{\vec{\alpha},q}{\Delta}A^{q,\vec{\alpha}}_{i,mc(X)}$$ Let $p^*\leq^* p'$ be the condition obtained by shrinking to those sets. $p'$ has the property that whenever $(*)_{q,X}$ holds for some $q\in\mathbb{Q}$ and $X\in Ex(p')$,
there exists sets $A(q,\vec{\alpha})$ for $\vec{\alpha}\in (X\setminus\{mc(X)\})(p')$ such that for every $ \vec{\alpha}^{\smallfrown}\l\alpha\r\in X(p')$, $A(q,\vec{\alpha})\cap\alpha=a(q,\vec{\alpha},\alpha)$.  By density there is such $p'\in G_2$.

Work in $V[A]$, for every $\vec{\alpha}$ and $q$, if $A(q,\vec{\alpha})$ is defined, let $$\eta(q,\vec{\alpha})=\min(A\Delta A(q,\vec{\alpha}))$$ otherwise $\eta(q,\vec{\alpha})=0$. $\eta(q,\vec{\alpha})$ is well defined since $A\notin V[C^*]$ and $A(q,\vec{\alpha})\in V[C^*]$. Also let $$\eta(\vec{\alpha})=\sup(\eta(q,\vec{\alpha})\mid q\in\mathbb{Q})$$ If $\eta(\vec{\alpha})=\kappa$ then we are done (since $|\mathbb{Q}|<\kappa$).  Define a sequence in $V[A]$:
$\alpha_0=\kappa^*$.  Fix $\xi<\otp(C_G)$ and assume that $\langle\alpha_i\mid i<\xi\rangle$ is defined. At limit stages take $$\alpha_\xi=\sup(\alpha_i\mid i<\xi)+1$$ Assume that $\xi=\lambda+1$ and let $$\alpha_\xi=\sup(\eta(\vec{\alpha})+1 \mid\vec{\alpha}\in[\alpha_\lambda]^{<\omega})$$ If at some point we reach $\kappa$ we are done. If not, let us prove by induction on $\xi$ that $C_G(\xi)<\alpha_\xi$ which will indicate that the sequence $\alpha_\xi$ is unbounded in $\kappa$. At limit $\xi$ we have $C_G(\xi)=\sup(C_G(\beta)\mid \beta<\xi)$ since the Magidor sequence is a club. By the definition of the sequence $\alpha_\xi$ and the induction hypothesis, $\alpha_\xi>C_G(\xi)$. If $\xi=\lambda+1$, use corollary \ref{IndCG} to find $\vec{\alpha}{}^{\smallfrown}\l\alpha\r$ and $q\in\mathbb{Q}$ such that 
$$\langle q,p'{}^{\frown}\vec{\alpha}{}^{\smallfrown}\l\alpha\rangle\rangle\Vdash\check{\alpha}=\underaccent{\sim}{C}_G(\check{\xi})$$
Fix any $q'\in \mathbb{Q}$ above $q$, and split the forcing at $\alpha$ so that $\langle q',p'^{\frown}\vec{\alpha}{}^{\smallfrown}\l\alpha\r\rangle=\langle q',r_1,r_2\rangle$ where $r_1\in\Mfor\restriction(\kappa^*,\alpha)$ and $r_2\in\Mfor\restriction(\alpha,\kappa)$. Let $H_1$ be some generic up to $\alpha$ with $\langle q',r_1\rangle\in H_1$ and work in $V[C^*][H_1]$, the name $\underaccent{\sim}{A}$ has a natural interpretation in $V[C^*][H_1]$ as a $\Mfor\restriction(\alpha,\kappa)$-name, $(\underaccent{\sim}{A})_{H_1}$. Use the fact that $\Mfor\restriction(\alpha,\kappa)$ is $\leq^*$-closed and the Prikry condition to find $r_2\leq^* r_2'\in \Mfor\restriction(\alpha,\kappa)$ and $A_0$ such that $$r_2'\Vdash_{\Mfor\restriction(\alpha,\kappa)} (\underaccent{\sim}{A})_{H_1}\cap\alpha=A_0$$ since it is forced that $\underaccent{\sim}{A}$ is old, $A_0\in V[C^*]$ and therefore we can find $\langle q'',r_1'\rangle\in \mathbb{Q}\times \Mfor\restriction(\kappa^*,\alpha)$ such that $\langle q'',r_1'\rangle\geq\langle q',r_1\rangle$ and
$$\langle q'',r_1'\rangle\Vdash ``r_2'\Vdash \underaccent{\sim}{A}\cap \alpha=A_0"\text{ Therefore } \langle q'',r_1',r_2'\rangle\Vdash \underaccent{\sim}{A}\cap\alpha=A_0$$

Since $r_2\leq^* r_2'$ and $r'_1\in\Mfor\restriction(\kappa^*,\alpha)$, then there is some $\vec{\beta}\in[\alpha]^{<\omega}$ such that $\langle r_1',r_2'\rangle{}^*\geq p'^{\frown}\vec{\beta}{}^{\smallfrown}\l\alpha\r $. Let $X$ be the extension type of $\vec{\beta}{}^{\smallfrown}\l\alpha\r$,  by definition of $p'$, $(*)_{q'',X}$ holds. Use density to find a condition $q^*$ in the generic of $\mathbb{Q}$ 
such that for some extension type $X$ that decides the $\xi$th element of $C_G$, $(*)_{X,q^*}$ holds. The set $\{p'^{\frown}\vec{\gamma}\mid \vec{\gamma}\in X(p')\}$ is a maximal antichain according to
proposition \ref{Max}, so let $\vec{C}
{}^{\smallfrown}C_G(\xi)$ be the extension of $p'$ of type $X$ in $C_G$. By the construction of $q^*$ and $p'$ we have that $$\langle q^*,p'{}^{\smallfrown} \vec{C}{}^{\smallfrown}C_G(\xi)\rangle \Vdash \ \underaccent{\sim}{A}\cap\check{C_G(\xi)}=A(q^*,\vec{C})\cap\check{C_G(\xi)}$$ Since $(\underaccent{\sim}{A})_G=A$, $A(q^*,\vec{C})\cap C_G(\xi)=A\cap C_G(\xi)$ (otherwise we would've found compatible conditions forcing contradictory information). This implies that $$\eta(q^*,\vec{C})\geq C_G(\xi)$$ By the induction hypothesis $\alpha_\lambda>C_G(\lambda)$ and $\vec{C}\subseteq C_G(\lambda)$ thus $\vec{C}\in [\alpha_{\lambda}]^{<\omega}$ thus
 $$\alpha_\xi>sup(\eta(\vec{\alpha})\mid \vec{\alpha}\in [\alpha_{\lambda}]^{<\omega})\geq\eta(\vec{C})\geq \eta(q^*,\vec{C})\geq C_G(\xi)$$
This proves that $\langle\alpha_{\xi}\mid \xi<\otp(C_G)<\kappa\rangle\in V[A]$ is cofinal in $\kappa$ indicating  $cf^{V[A]}(\kappa)<\kappa$.
$\blacksquare$

Thus we have proven the result for any subset of $\kappa$.
\begin{corollary}\label{Carkap}
Let $A\in V[G]$ be a set of ordinals, such that $|A|=\kappa$ then there is $C'\subseteq C_G$ such that $V[A]=V[C']$. 
\end{corollary}
\pr By $\kappa^+$-c.c. of $\Mfor$, there is $B\in V$, $|B|=\kappa$ such that $A\subseteq B$. Fix in $V$ $\phi:\kappa\rightarrow B$ a bijection and let $B'=\phi^{-1''}A$. then $B'\subseteq\kappa$. By the theorem for subsets of $\kappa$ there is $C'\subseteq C_G$ such that $V[C']=V[B']=V[A]$.
$\blacksquare$
\subsection{General sets of ordinals}

  In \cite{TomMoti}, we gave an explicit formulation of subforcings of $\Mfor$ using the indices of subsequences of $C_G$. In the larger framework of this paper, these indices might not be in $V$. By example \ref{examplenonGeneralize},  subforcing of the Magidor forcing can be an iteration of Magidor type forcing. 
\begin{lemma}
Let $A\in V[G]$ be such that $A\subseteq\kappa^+$. Then there is $C^*\subseteq C_G$ such that
\begin{enumerate}
\item $\exists\alpha^*<\kappa^+$ such that $C^*\in V[A\cap\alpha^*]\subseteq V[A]$.
\item $\forall\alpha<\kappa^+ \ A\cap\alpha\in V[C^*]$.
\end{enumerate}
\end{lemma}
\pr Work in $V[G]$, for every $\alpha<\kappa^+$ find subsequences $C_\alpha\subseteq C_G$  such that $$V[C_\alpha]=V[A\cap\alpha]$$ using corollary \ref{Carkap}.
The function $\alpha\mapsto C_\alpha$ has range $P(C_G)$ and domain $\kappa^+$ which is regular in $V[G]$, and since $o^{\vec{U}}(\kappa)<\kappa$ then $|C_G|<\kappa$, and since $\kappa$ is strong limit (even in $V[G]$) $|P(C_G)|<\kappa<\kappa^+$.
Therefore there exist $E\subseteq\kappa^+$ unbounded in $\kappa^+$ and $\alpha^*<\kappa^+$  such that for every $\alpha\in E$, $C_\alpha=C_{\alpha^*}$. Set $C^*=C_{\alpha^*}$, note that for every $\alpha<\kappa$ there is $\beta\in E$ such that $\beta>\alpha$ therefore $$A\cap\alpha=(A\cap\beta)\cap\alpha\in V[A\cap\beta]=V[C^*]$$
$\blacksquare$

\begin{lemma}\label{suffcondidi}
Let $C^*$ be as in the last lemma. If there is  $\alpha<\kappa$ such that $A\in V[C_G\cap\alpha][C^*]$ then $V[A]=V[C^*]$. \end{lemma}
\pr Consider the quotient forcing $\Mfor/C^*\subseteq \Mfor$ completing $V[C^*]$ to $V[C^*][G]$. Then the forcing  $$\mathbb{Q}=(\Mfor/C^*)\restriction\alpha$$ completes $V[C^*]$ to $V[C^*][C_G\cap\alpha]$ and $|\mathbb{Q}|<\kappa$. By the assumption, $A\in V[C^*][C_G\cap\alpha]$, and for every $\beta<\kappa^+$, $A\cap\beta\in V[C^*]$. Let $\underaccent{\sim}{A}\in V[C^*]$ be a $\mathbb{Q}$-name for $A$ and $q\in G\restriction \alpha$ be any condition such that
$$q\Vdash \forall\beta<\kappa^+,\underaccent{\sim}{A}\cap\beta \in V[C^*]$$
In $V[C^*]$, for every $\beta<\kappa^+$ find $q_\beta\geq q$ such that $q_\beta ||_{\mathbb{Q}} \underaccent{\sim}{A}\cap\beta$, there is $q^*\geq q$ and $E\subseteq\kappa^+$ of cardinality $\kappa^+$ such that for very $\beta\in E$, $q_\beta=q^*$. 
By density, find such $q^*\in G\restriction \alpha$ in the generic. In $V[C^*]$, consider the set $$B=\{X\subseteq\kappa^+\mid \exists\beta \ q^*\Vdash X=\underaccent{\sim}{A}\cap\beta\}$$
Let us argue that $\cup B=A$. Let $X\in B$  then there is $\beta<
\kappa^+$ such that $q^*\Vdash X=\underaccent{\sim}{A}\cap \beta$ then $X=A\cap\beta\subseteq A$, thus, $\cup B\subseteq A$. Let $\gamma\in A$, there is $\beta\in E$ such that $\gamma<\beta$, by the definition of $E$ there is $X\subseteq\beta$ such that $q^*\Vdash\underaccent{\sim}{A}\cap\beta=X$ it must be that $X=A\cap\beta$ otherwise would have found compatible conditions forcing contradictory information. But then $\gamma\in A\cap\beta=X\subseteq\cup B$. We conclude that $A=\cup B\in V[C^*]$.
$\blacksquare$
\vskip 0.3 cm
Eventually we will prove that there is $\alpha<\kappa$ such that $A\in V[C_G\cap\alpha][C^*]$ and by the last lemma we will be done.

We would like to change $C^*$ so that it is closed. We can do that above $\alpha_0:=\otp(C_G)$:
\begin{lemma}\label{Clubset}
 $V[C_G\cap\alpha_0][Cl(C^*)]=V[C_G\cap\alpha_0][C^*]$.\footnote{ For a set of ordinals $X$, $Cl(X)=X\cup \Lim(X)=\{\xi\mid \xi\in X\vee \sup(X\cap\xi)=\xi\}   $}
\end{lemma}
\pr 
Consider $I(C^*, Cl(C^*))\subseteq \otp(C_G)$, by proposition \ref{genericproperties}(5), $I(C^*,Cl(C^*))\in V[C_G\cap\alpha_0]$.
Thus $V[C_G\cap\alpha_0][C^*]=V[C_G\cap\alpha_0][Cl(C^*)]$.
$\blacksquare$

Work in $N:=V[C_G\cap\alpha_0]$, since $C^*\cap\alpha_0\in V[C_G\cap\alpha_0]$, we can assume $min(C^*)>\alpha_0$. Since $I=I(C^*,C_G\setminus\alpha_0)\subseteq \otp(C_G)$, it follows that $I\in N$. In $N$, consider the coherent sequence $$\vec{W}=\vec{U}^*\restriction(\alpha_0,\kappa]=\langle U^*(\beta,\delta)\mid\delta<o^{\vec{U}}(\beta), \alpha_0<\delta<\kappa\rangle$$ where $U^*(\beta,\delta)$ is the ultrafilter generated by $U(\beta,\delta)$ in $N$. Also denote $G^*=G\restriction(\alpha_0,\kappa)$. The following proposition is to be compared with remark \ref{lift}.

\begin{proposition}
 $N[G^*]$ is a $\mathbb{M}[\vec{W}]$ generic extension of $N$.
\end{proposition}
\pr Let us argue that the Mathias criteria holds. Let $X\in\cap\vec{W}(\delta)$ where $\delta\in Lim(C_{G^*})$. By definition of $\vec{W}$, for every $i<o^{\vec{W}}(\delta)$, there is $X_i\in U(\delta,i)$, such that $X_i\subseteq X$. The choice of $X_i$'s is done in $N$ and the sequence $\l X_i\mid i<o^{\vec{U}}(\delta)\r$ might not be in $V$. Fortunately, $\Mfor\restriction \alpha_0$ is $\alpha_0^+$-c.c. and $\alpha_0^+<\delta$, so in $V$, we can find sets $$E_i:=\{X_{i,j}\mid j\leq \alpha_0\}\subseteq U(\delta,i)$$ such that 
$X_i\in E_i$ By $\delta$-completness of $U(\delta,i)$, the set $X_i^*:=\cap E_i\in U(\delta,i)$ and $X_i^*\subseteq X_i\subseteq X$. Note that $X^*:=\cup_{i<o^{\vec{U}}(\delta)}X_i^*\in\cap \vec{U}(\delta)$ and therefore by genericity of $G$ there is $\xi<\delta$ such that $$C_G\cap (\xi,\delta)\subseteq X^*\subseteq X$$
Hence $C_{G^*}\cap(\max(\alpha_0,\xi),\delta)\subseteq X$.$\blacksquare$

Note that $o^{\vec{W}}(\kappa)<\min\{\nu\mid o^{\vec{W}}(\nu)=1\}$ and $I(C^*,C_G)\in N$, In \cite{TomMoti}, this is the situation dealt with, a forcing denoted by $\MforIW\in N[C^*]$ was defined where $I=I(C^*,C_G)$ and used to conclude the theorem. We only state here the main results and definitions and refer the reader to \cite{TomMoti} for the full definition and proofs. 

\begin{proposition}\label{Whatwe}
Let $G^*\subseteq \mathbb{M}[\vec{W}]$ be $N$-generic filter and $C\subseteq C_{G^*}$ be closed. Assume that $I=I(C,C_{G^*})\in N$. Then there is a forcing notion $\MforIW\in N$ and a projection $\pi_I:\mathbb{M}[\vec{W}]\rightarrow\MforIW$ such that $N[G_I]=N[C]$ where $G_I=\overline{\pi_I''G^*}\subseteq\MforIW$ is the $N$-generic filter obtained by projecting $G^*$.
\end{proposition}
\begin{lemma}\label{ccG}
Let $G^*\subseteq \mathbb{M}[\vec{W}]$ be $N$-generic filter. Then the forcing $\mathbb{M}[\vec{W}]/G_I$ satisfies $\kappa^+-c.c.$ in $N[G^*]$.
\end{lemma}
The referee pointed out a simpler argument then the one given in \cite{TomMoti} for the continuation of the proof.
First we conclude the following (see for example \cite[Thm 16.4]{Jech2003}:
\begin{corollary}\label{prodisc}
The forcing $\mathbb{M}[\vec{W}]/G_I\times \mathbb{M}[\vec{W}]/G_I$ satisfy $\kappa^+$-c.c.
\end{corollary}
The next theorem is what needed in order to apply lemma \ref{suffcondidi} and to conclude the case for $A\subseteq\kappa^+$.
\begin{theorem}
$A\in N[C^*]$.
\end{theorem}
\pr
 Let $I=I(Cl(C^*),C_{G^*})$. Then $$I,\MforIW,\pi_I\in N$$ Let $G_I$ be the generic induced for $\MforIW$ from $G$, it follows that $\mathbb{M}[\vec{W}]/G_I$ is defined in $N$. Toward a contradiction, assume that $A\notin N[C^*]$. By lemma \ref{Clubset}, $N[C^*]=N[Cl(C^*)]$, hence $A\notin N[Cl(C^*)]$. Let $\underaccent{\sim}{A}$ be a name for $A$ in $\Mfor/G_I$. Work in $N[G_I]$, by corollary \ref{Whatwe}, $N[G_I]=N[Cl(C^*)]$.  We define a tree $T\in N[G_I]$ of height $\kappa^+$, For every $\alpha<\kappa^+$ define the $\alpha$th level of the tree by
$$Lev_\alpha(T)=\{B\subseteq\alpha \mid ||\underaccent{\sim}{A}\cap\alpha=B||\neq 0\}$$
where the truth value is taken in $RO(\mathbb{M}[\vec{W}]/G_I)$- the complete Boolean algebra of regular open sets for  $\mathbb{M}[\vec{W}]/G_I$. The order of the tree $T$ is simply end- extension. 
Different $B$'s in $Lev_\alpha(T)$ yield incompatible conditions of $\mathbb{M}[\vec{W}]/G_I$ and we have $\kappa^+$-c.c by lemma \ref{ccG} thus
$$\forall\alpha<\kappa^+ \ |Lev_\alpha(T)|\leq\kappa$$

 Work in $N[G^*]$,  denote $A_\alpha=A\cap \alpha$. Recall that
$$\forall \alpha<\kappa^+ \ A_\alpha\in N[Cl(C^*)]=N[G_I]$$
thus $A_\alpha\in Lev_\alpha(T)$ which makes $A$ a branch through $T$. At this point, the referee pointed out an argument by Unger \cite{Unger2013} showing that a forcing $\mathbb{P}$ such that $\mathbb{P}\times\mathbb{P}$ satisfy $\kappa^+$-c.c. has the $\kappa^+$-approximation property and in particular cannot add a new branches to $\kappa^+$ trees in the ground model (see definition $2.2$, the discussion succeeding it, and lemma $2.4$ in \cite{Unger2013}). By corollary \ref{prodisc}, the product of $\mathbb{M}[\vec{W}]/G_I$ in $\kappa^+-c.c$ in $N[G_I]$ and therefore $\mathbb{M}[\vec{W}]/G_I$ does not add new branches to $\kappa^+$ which implies that $A\in N[G_I]$.

For self inclusion reasons and for the convenience of the reader, let us give another argument, for every $B\in Lev_\alpha(T)$ define $$b(B)=||\underaccent{\sim}{A}\cap\alpha=B||$$
Assume that $B'\in Lev_\beta(T)$ and $\alpha\leq\beta$ then $B=B'\cap\alpha \in Lev_\alpha(T)$. Moreover $b(B')\leq_B b(B)$ (we switch to Boolean algebra notation $p\leq_B q$ means $p$ extends $q$). Note that for such $B,B'$ if $b(B')<_B b(B)$, then there is 
$$0<p\leq_B(b(B)\setminus b(B'))\leq_Bb(B)$$
Therefore
 $$p\cap b(B')\leq_B (b(B)\setminus b(B'))\cap b(B')=0$$
meaning $p\bot b(B')$. As before, in $N[G^*]$ we denote $A_\alpha=A\cap \alpha\in Lev_\alpha(T)$.
 Consider the  $\leq_B$-non-increasing sequence $\langle b(A_\alpha) \mid \alpha<\kappa^+\rangle$. If there exists some $\gamma^*<\kappa^+$ on which the sequence stabilizes, define 
$$A'=\bigcup\{B\subseteq\kappa^+ \ | \ \exists\alpha \ b(A_{\gamma^*})\Vdash \underaccent{\sim}{A}\cap\alpha=B\}\in N[Cl(C^*)]$$
Claim that $A'=A$, notice that if $B,B',\alpha,\alpha'$ are such that 
$$b(A_{\gamma^*})\Vdash \underaccent{\sim}{A}\cap\alpha=B, \ \  
b(A_{\gamma^*})\Vdash \underaccent{\sim}{A}\cap\alpha'=B'$$
With out loss of generality, $\alpha\leq\alpha'$ then we must have $B'\cap\alpha=B$ otherwise, the non zero condition $b(A_{\gamma^*})$ would force contradictory information. Consequently, for every
$\xi<\kappa^+$ there exists $\xi<\gamma<\kappa^+$ such that $$b(A_{\gamma^*})\Vdash \underaccent{\sim}{A}\cap\gamma=A\cap\gamma$$ hence $A'\cap\gamma= A\cap\gamma$. This is a contradiction to $A\notin N[Cl(C^*)]$. We conclude that he sequence $\langle b(A_\alpha) \mid \alpha<\kappa^+\rangle$ does not stabilize. By regularity of $\kappa^+$, there exists a subsequence $$\langle b(A_{i_\alpha}) \mid \alpha<\kappa^+\rangle$$ which is strictly decreasing. Use the observation we made to find $p_\alpha\leq_B b(A_{i_\alpha})$ such that $p_\alpha \bot b(A_{i_{\alpha+1}})$. Since $b(A_{i_\alpha})$ are decreasing, for any $\beta>\alpha \  p_\alpha\bot b(A_{i_\beta})$ thus $p_\alpha\bot p_\beta$. This shows that $\langle p_\alpha \mid \alpha<\kappa^+\rangle\in N[G^*]$ is an antichain of size $\kappa^+$ which contradicts Lemma \ref{ccG}.
$\blacksquare$

\vskip 0.3 cm
\textbf{Sets of ordinals above $\kappa^+$}: By induction on $\sup(A)=\lambda>\kappa^+$. It suffices to assume that $\lambda$ is a cardinal. 
\vskip 0.2 cm
\underline{case1:} $cf^{V[G]}(\lambda)>\kappa$, the arguments for $\kappa^+$ works.
\vskip 0.2 cm
\underline{case2:} $cf^{V[G]}(\lambda)\leq\kappa$ and since $\kappa$ is singular in $V[G]$ then $cf^{V[G]}(\lambda)<\kappa$. Since $\Mfor$ satisfies $\kappa^+-c.c.$ we must have that   $\nu:=cf^V(\lambda)\leq\kappa$. Fix $$\langle\gamma_i | \ i<\nu\rangle\in V$$ cofinal in $\lambda$. Work in $V[A]$, for every $i<\nu$ find $d_i\subseteq \kappa$ such that $V[d_i]=V[A\cap\gamma_i]$. By induction, there exists $C^*\subseteq C_G$ such that $V[\langle d_i\mid i<\nu\rangle]=V[C^*]$, therefore 
\begin{enumerate}
\item $\forall i<\nu \ A\cap\gamma_i\in V[C^*]$
\item $C^*\in V[A]$
\end{enumerate}
Work in $V[C^*]$, for $i<\nu$ fix $$\langle X_{i,\delta}\mid \delta<2^{\gamma_i}\rangle=P(\gamma_i)$$ then we can code $A\cap\gamma_i$ by some $\delta_i$ such that $X_{i,\delta_i}=A\cap\gamma_i$. By \ref{Carkap}, we can find $C''\subseteq C_G$ such that $$V[C'']=V[\langle \delta_i\mid i<\nu\rangle]$$ Finally we can find $C'\subseteq C_G$ such that $V[C']=V[C^*,C'']$, it follows that
$V[A]=V[C']$.
$\blacksquare_{\text{Theorem }\ref{MainResaultPartone}}$
\section{Classification of Intermediate Models}

Let $G\subseteq\Mfor$ be a $V$-generic filter. Assume that for every $\alpha\leq\kappa$, $o^{\vec{U}}(\alpha)<\alpha$. Let $M$ be a transitive $ZFC$ model such that $V\subseteq M\subseteq V[G]$.
We would like to prove it is a generic extension of a ``Magidor-like" forcing which will be defined shortly.

By example \ref{examplenonGeneralize}, the class of forcings $\MforI$ does not capture all the intermediate models of a generic extension by $\Mfor$. The reason is that if $o^{\vec{U}}(\kappa)\geq \min\{\alpha\mid o^{\vec{U}}(\alpha)=1\}$, there are subsets $C\subseteq C_G$ such that $I(C,C_G)$ does not necessarily exists in the ground model, which was crucial in the definition of $\MforI$. Here we generalize this class to a class of forcings denoted by $\mathbb{M}_f[\vec{U}]$. We will prove that every intermediate model is a generic extension for a finite iteration of forcings of the form $\mathbb{M}_f[\vec{U}]$. The major difference between $\mathbb{M}_f[\vec{U}]$ and $\MforI$ is the existence of a concrete projection of $\Mfor$ onto $\MforI$  which keeps only the ordinals which will sit at index $i\in I$ in the generic club. As for the generic set produced by $\mathbb{M}_f[\vec{U}]$, we cannot determine in advance how this set sits inside $C_G$. For example if $\MforI$ turns out to be the standard Prikry forcing, then the projection tells us what indices the Prikry sequence fill in $C_G$, and the forcing made sure to leave ``room" for the missing elements of $C_G$. On the other hand, if  $\mathbb{M}_f[\vec{U}]$ produces a Prikry sequence, there will be many ways to place this Prikry sequence inside $C_G$. One might claim that this is only a technicality, but if we aim to describe a forcing which produces a generic extension for an intermediate model of the form $V[C]$, where $C\subseteq C_G$, then example \ref{IndexNotInVC} below describe a situation that $I(C,C_G)\notin V[C]$, and in particular there is no model $V\subseteq N\subseteq V[C]$ such that $V[C]$ is a generic extension of $N$ by $\MforI$. Instead of using $I(C,C_G)$, the forcing $\mathbb{M}_f[\vec{U}]$ uses the sequence $\l o^{\vec{U}}(\alpha) \mid \alpha\in C\r$ which is definable in $V[C]$.

\begin{example}\label{IndexNotInVC}
Consider $\kappa$ such that $o^{\vec{U}}(\kappa)=\delta_0:=\min\{\alpha\mid o^{\vec{U}}(\alpha)=1\}$. Let $$p=\l \l\delta_0,A\r,\l\kappa,B\r\r\in\Mfor$$
then $p\Vdash C_{\lusim{G}}(\omega)=\delta_0$.
Let $G\subseteq \Mfor$ be such that $p\in G$, and consider the first Prikry sequence for $C_G(\omega)=\delta_0$, namely $\{C_G(n)\mid n<\omega\}$, and let $C=\{C_G(C_G(n)+1)\mid n<\omega\}$. Since for each $n<\omega$, $C_G(C_G(n)+1)$ is successor in $C_G$,  $o^{\vec{U}}(C_G(C_G(n)+1))=0$ and therefore $C$ is a Prikry sequence for $U(\kappa,0)$.
Note that $I(C,C_G)=\{C_G(n)+1\mid n<\omega\}$ and $I(C,C_G)\notin V[C]$. Otherwise $\{C_G(n)\mid n<\omega\}\in V[C]$, which is a contradiction since Prikry extensions do not add bounded subsets to $\kappa$.
\end{example}

\begin{proposition}\label{union of generic}
Let $C,D\subseteq C_G$. There exists $E$ such that $C\cup D\subseteq E\subseteq C_G\cap \sup( C\cup D)$ and $V[C,D]=V[E]$.
\end{proposition}
\pr By induction on $\sup(C\cup D)$. If $\sup(C\cup D)\leq C_G(\omega)$ then $|C|,|D|\leq\aleph_0$, we can take $E=C\cup D$, clearly $$I(C,C\cup D),I(D,C\cup D)\subseteq\omega$$ thus these sets belong to $V$.
In the general case, consider $I(C,C\cup D),I(D,C\cup D)$. Since $$o^{\vec{U}}(\sup(C\cup D))<\sup(C\cup D)$$
it follows that
$$\otp(C\cup D)\leq \otp(C_G\cap\sup(C\cup D))<\sup(C\cup D)$$ Denote by $\lambda=\otp(C_G\cap\sup(C\cup D))$. By theorem \ref{MainResaultPartone}, there is $F\subseteq C_G\cap \lambda$
,  such that $$V[I(C,C\cup D),I(D,C\cup D)]=V[F]$$ 
Apply the induction hypothesis to $F, (C\cup D)\cap \lambda$ and find $E_*\subseteq \lambda$ such that $$V[E_*]=V[F, (C\cup D)\cap \lambda]$$ Let $E=E_*\cup (D\cup C)\setminus \lambda$, then $E\in V[C,D]$ as both $E_*,D\cup C$ in $V[C,D]$. In $V[E]$ we can find $$E_*=E\cap \lambda\text{ and }(D\cup C)\setminus \lambda=E\setminus \lambda$$ Thus $F,(C\cup D)\cap \lambda\in V[E]$ and therefore also  $$D\cup C, I(C,C\cup D),I(D,C\cup D)\in V[E]$$ It follows that $C,D\in V[E]$.$\blacksquare$
\begin{corollary}\label{FindingClosedSet}
For every $C'\subseteq C_G$ there is $C^*\subseteq C_G\cap\sup(C')$, such that $C^*$ is closed and $V[C']=V[C^*]$.
\end{corollary}
\pr Again we go by induction on $\sup(C')$.
If $\sup(C')=C_G(\omega)$ then $C^*=C'$ is already closed.
For general $C'$, consider $C'\subseteq Cl(C')$, then $I(C',Cl(C'))$ is bounded by some $\nu<\sup(C')$. So there is $D\subseteq C_G\cap\nu$ such that $V[D]=V[I(C',Cl(C'))]$. By proposition \ref{union of generic}, we can find $E$ such that $$D\cup Cl(C')\cap\nu\subseteq E\subseteq C_G\cap\nu$$ and $V[E]=V[D,Cl(C')]$. By the induction hypothesis there is a closed $E_*$, such that $E\subseteq E^*\subseteq C_G\cap \nu$ and $V[E]=V[E_*]$.
Finally, let $$C^*=E_*\cup\{\sup(E_*)\}\cup Cl(C')\setminus\nu$$
Then $C^*\in V[C']$, and also $Cl(C')$ and $I(C',Cl(C'))$ can be constructed in $V[C^*]$ so $C'\in V[C^*]$. Obviously, $C^*$ is closed, hence, $C^*$ is as desired.$\blacksquare$
\begin{definition}
Let $\lambda<\kappa$ be ordinal. A function $f:\lambda\rightarrow \kappa$ is \textit{suitable}, if for all $\delta\in Lim(\lambda)$,
$$\limsup_{\alpha<\delta}f(\alpha)+1\leq f(\delta)$$
\end{definition}

We would like to define $\mathbb{M}_f[\vec{U}]$ for a suitable $f$, to be the forcing which construct a continuous sequence such that the order of the elements of the sequence is prescribed by $f$. However we must require some connection to $\vec{U}$. In example \ref{someconnection} below, we provide a suitable function which cannot describe the orders of any generic subsequence. 
\begin{example}\label{someconnection}
Assume that $o^{\vec{U}}(\kappa)=\omega_1$ and  $\forall\alpha<\kappa.o^{\vec{U}}(\alpha)<\omega_1$. 
Let $f:\omega+1\rightarrow \kappa$ defined by $f(0)=f(\omega)=\omega_1$ and $f(n+1)=0$.
 There is no $C\subseteq C_G\cup\{\kappa\}$ with $\otp(C)=\omega+1$ such that $o^{\vec{U}}(C(i))=f(i)$. There are two reasons for that, the first, is that there is no $\alpha<\kappa$ that can be $C(0)$, since by assumption $o^{\vec{U}}(\alpha)<\omega_1=f(0)$. The second reason is that $cf^{V[G]}(\kappa)=\omega_1$, hence there is no unbounded $\omega$-sequence of ordinals of order $0$ below $\kappa$.
\end{example}

Let us restrict our attention to a more specific family of suitable functions.
\begin{definition}
Let $G\subseteq\Mfor$ be $V$-generic and $C\subseteq C_G$ be closed,  $\lambda+1=\otp(C\cup\{\sup(C)\})$, and $\langle C(i)\mid i\leq\lambda\rangle$ be the increasing continuous enumeration of $C$. \textit{The suitable function derived from $C$} denoted by $f_{C}$, is the function 
$f_C:\lambda+1\rightarrow \kappa$, defined by $f_C(i)=o^{\vec{U}}(C(i))$.
A suitable function is called a \textit{derived suitable function} if it is derived from some closed $C\subseteq C_G$.
\end{definition}

\begin{proposition}
If $C\subseteq C_G$ is a closed subset,  then $f_C$  is suitable.
\end{proposition}
\pr Let $\delta\in\Lim(\lambda+1)$, then $C(\delta)\in Lim(C_G\cup\{\kappa\})$ and therefore, there is $\xi<C(\delta)$ such that for every $x\in C_G\cap(\xi,C(\delta))$, $o^{\vec{U}}(x)<o^{\vec{U}}(C(\delta))$. Let $\rho<\delta$ be such that for every $\rho<i<\delta$, $\xi<C(i)<C(\delta)$. Then $\sup_{\rho<i<\delta}o^{\vec{U}}(C(i))+1\leq o^{\vec{U}}(C(\delta))$ and also $$\min\Big\{(\sup_{\alpha<i<\delta}o^{\vec{U}}(C(i))+1)\mid \alpha<\delta\Big\}\leq o^{\vec{U}}(C(\delta))$$ $\blacksquare$

\begin{definition}\label{DefinitionMfU}
Let $f:\lambda+1\rightarrow \kappa$ be a derived suitable function. Define the forcing $\mathbb{M}_f[\vec{U}]$, the conditions are functions $F$, such that:
\begin{enumerate}
    \item  $F$ is finite partial function, with $Dom(F)\subseteq \lambda+1$. such that $\lambda\in Dom(F)$.
    \item For every $ i\in Dom(F)\cap Lim(\lambda+1)$:
    \begin{enumerate}
        \item $F(i)=\langle \kappa^{(F)}_i,A^{(F)}_i\rangle$.
        \item $o^{\vec{U}}(\kappa^{(F)}_i)=f(i)$.
        \item $A^{(F)}_i\in\cap\vec{U}(\kappa^{(F)}_i)$.
        \item Let $j=\max(Dom(F)\cap i)$ or $j=-1$ if $i=\min(Dom(F))$, then for every $j<k<i$, $f(k)<f(i)$.
    \end{enumerate}
    \item For every $i\in Dom(F)\setminus Lim(\lambda)$
    \begin{enumerate}
        \item $F(i)=\kappa^{(F)}_i$.
        \item $o^{\vec{U}}(\kappa^{(F)}_i)=f(i)$.
        \item $i-1\in Dom(F)$. 
    \end{enumerate}
    \item The map $i\mapsto \kappa^{(F)}_i$ is increasing.
    
\end{enumerate} 
\end{definition}
\begin{definition}
The order of $\mathbb{M}_f[\vec{U}]$ is defined as follows $F\leq G$ iff
\begin{enumerate}
    \item $Dom(F)\subseteq Dom(G)$.
    \item For every $i\in Dom(G)$, let $j=\min(Dom(F)\setminus i)$. 
    \begin{enumerate}
        \item If $i\in Dom(F)$, then $\kappa^{(F)}_i=\kappa^{(G)}_i$, and $A^{(G)}_i\subseteq A^{(F)}_i$.
    \item If $i\notin Dom(F)$, then $\kappa^{(G)}_i\in A^{(F)}_j$, and $A^{(G)}_i\subseteq A^{(F)}_j$.
    \end{enumerate}
\end{enumerate}
\end{definition}
\begin{proposition}\label{Forcingnotion}
 Let $f$ be a suitable derived function, then $\mathbb{M}_f[\vec{U}]$ is a forcing notion.
\end{proposition}
\pr It is not hard to check that $\leq$ is a partial order on $\mathbb{M}_f[\vec{U}]$. To see $\mathbb{M}_f[\vec{U}]\neq\emptyset$, let $C$ be such that $f=f_C$. We define a finite sequence $\alpha_0=\lambda$,if $\alpha_0$ is successor, $\alpha_1=\alpha_0-1$. Otherwise, if there is no $\beta$ such that $f(\beta)\geq f(\alpha_0)$, then we halt the definition. If there is such $\beta$, let $\alpha_1=\max\{\beta<\alpha_0\mid f(\beta)\geq f(\alpha_0)\}$. By the suitability requirement, this maximum is defined and $\alpha_1<\alpha_0$. In a similar fashion if $\alpha_1$ is successor, let $\alpha_2=\alpha_1-1$, if there is no $\beta$ such that $f(\beta)\geq f(\alpha_1)$, then we halt the definition, otherwise, $\alpha_2=\max\{\beta<\alpha_1\mid f(\beta)\geq f(\alpha_1)\}$ and $\alpha_2<\alpha_1<\alpha_0$. After finitely many steps we reach $\alpha_k$ such that for every $\beta<\alpha_k$, $f(\beta)<f(\alpha_k)$.  The function $F$ defined by $Dom(F)=\{\alpha_k,...,\alpha_1\}$ and $F(\alpha_i)=\l C(\alpha_i),C(\alpha_i)\setminus C(\alpha_{i+1})+1\r$, satisfy definition \ref{DefinitionMfU}.$\blacksquare$

\begin{example}
Assume that $f:\omega+1\rightarrow \kappa$, defined by $f(n)=0$ and $f(\omega)=1$ then $\mathbb{M}_f[\vec{U}]$ first picks some measurable $\kappa^F_{\omega}$ of order $1$, then adds a Prikry sequence to the measure $U(\kappa^F_{\omega},0)$. 

If we only change $f$ at $\omega$, $f(\omega)=2$, then we still force a Prikry sequence for the measure $U(\kappa^F_{\omega},0)$, but the first part chooses a measurable of order $2$.
\end{example}
\begin{example}
Let $f:\omega^2+\omega+1\rightarrow \kappa$ defined by $f(\omega\cdot n+m)=n$, $f(\omega^2)=\omega$, $f(\omega^2+m+1)=1$ and $f(\omega^2+\omega)=2$. Clearly, $f$ is suitable. $\mathbb{M}_f[\vec{U}]$ first picks a measurable $\kappa^{(F)}_{\omega^2+\omega}$ of order $1$. By condition $2.d$ of definition \ref{DefinitionMfU},  we must also pick $\kappa^{(F)}_{\omega^2}$ of order $\omega$, since $f(\omega^2)>f(\omega^2+\omega)$. Then in the interval $(\kappa^{(F)}_{\omega^2},\kappa^{(F)}_{\omega^2+\omega})$ the forcing generates a Prikry sequence for $U(\kappa^{(F)}_{\omega^2+\omega},1)$ and below $\kappa^{(F)}_{\omega^2}$ the forcing generates an diagonal Prikry sequence $\{\kappa^{(F)}_{\omega^n}\mid n<\omega\}$ for the measures $\l U(\kappa^{(F)}_{\omega\cdot n},n)\mid n<\omega\r$. For each $n<\omega$, the forcing generates a Prikry sequence $\{\kappa^{(F)}_{\omega\cdot n+m}| \mid m<\omega\}$ for $U(\kappa^{(F)}_{\omega\cdot (n+1)},n)$ in the interval $[\kappa^{(F)}_{\omega\cdot n},\kappa^{(F)}_{\omega\cdot (n+1)})$. So in all $\mathbb{M}_f[\vec{U}]$ generates a sequence of order type $\omega^2+\omega+1$.
\end{example}

Let $f:\omega^{o^{\vec{U}}(\kappa)}+1\rightarrow \kappa$, defined by $f(\alpha)=o_L(\alpha)$ (see definition \ref{DefLimitOrder}). By proposition \ref{IndCG}, for every $V$-generic filter $G\subseteq\Mfor$ with $p_0:\langle\l\kappa,\kappa\r\r\in G$, $f=f_{C_G}$. Hence above $p_0$, $\Mfor$ is isomorphic to $\mathbb{M}_f[\vec{U}]$. Note that forcing with $\Mfor$ above $p_0$ is in the framework of this section since $\forall\alpha\in C_G\cup\{\kappa\}. o^{\vec{U}}(\alpha)<\alpha$.

Similar to $\Mfor$, we decompose sets $A^{(F)}_i=\biguplus_{\xi<o^{\vec{U}}(\kappa^{(F)}_i)}A^{(F)}_{i,\xi}$. Also if $j$ is as in condition $2.d$ of definition \ref{DefinitionMfU} and $j<i_1<...<i_k<i$, then for every $\vec{\alpha}\in\prod_{r=1}^kA^{(F)}_{f(i_r)}$, $G:=F^{\smallfrown}\vec{\alpha}$ is such that $Dom(G)=Dom(F)\cup\{i_1,...,i_k\}$ and $G(x)=F(x)$ unless $x=i_r$, in which case $G(x)=\vec{\alpha}(r)$.
\begin{proposition}
Let $f:\lambda+1\rightarrow \kappa\in V$ be a derived suitable function and $H\subseteq\mathbb{M}_f[\vec{U}]$ be a $V$-generic filter. Let $$C^*_H:=\{\kappa^{(F)}_i\mid i\in Dom(F),  F\in H\}$$
Then,
\begin{enumerate}
    \item $\otp(C^*_H)=\lambda+1$ and $C^*_H$ is continuous.
    \item For every $i\leq\lambda$, $o^{\vec{U}}(C^*_H(i))=f(i)$.
    \item $V[C^*_H]=V[H]$.
    \item For every $\delta\in\Lim(\lambda+1)$, and every $A\in\cap\vec{U}(\delta)$, there is $\xi<\delta$ such that $C^*\cap(\xi,\delta)\subseteq A$.
    \item For every successor $\rho<\lambda$, $H\restriction\rho:=\{ F\restriction\rho\mid F\in H\}$ is $V$-generic for $\mathbb{M}_{f\restriction\rho}[\vec{U}]$.
\end{enumerate}
\end{proposition}
\pr
To see $(1)$, let us argue  by induction on $i<\lambda$ that the set $$E_i=\{ F\in\mathbb{M}_{f}[\vec{U}]\mid i\in Dom(F)\}$$ is dense. Let $F\in \mathbb{M}_{f}[\vec{U}]$, if $i\in Dom(F)$ we are done. Otherwise,  let $$j_M:=\min(Dom(F)\setminus i)>i>\max(Dom(F)\cap i)=:j_m$$
By condition $3.c$ of definition \ref{DefinitionMfU} and minimality of $j_M$, $j_M\in Lim(\lambda+1)$. Split into two cases. First, if $i$ is successor, then we can find $F\leq G$ such that $i-1\in Dom(G)$ by induction hypothesis.
By condition $2.d$ and $2.b$, $f(i)<o^{\vec{U}}(\kappa_{j_M}^{(F)})$. By condition $2.c$, we can find $\alpha\in A^{(F)}_{j_M}$ such that $\alpha>\kappa^i_{j_m}$, $o^{\vec{U}}(\alpha)=f(i)$ and $A^{(F)}_{j_M}\cap\alpha\in\cap\vec{U}(\alpha)$. Then $$G'=G\cup\{\langle i,\langle \alpha, A^{(F)}_{j_M}\cap\alpha\rangle\r\}$$
is as wanted. If $i$ is limit, since $f$ is suitable, there is $i'<i$, such that for every $i'<k<i$, $f(k)<f(i)$. Again by induction, find $F\leq G$ such that $i'\in Dom(G)$. Then the desired $G'$ is constructed as in successor step. Denote by $F_H$, the function with domain $\lambda+1$, and $F_H(i)=\gamma$, be the unique $\gamma$ such that for some $ F\in H$, $i\in Dom(F)$ and\ $\kappa^{(F)}_i=\gamma$. Then it is clear that $F_H$ is order preserving and $1-1$ from $\lambda$ to $C^*_H$.
By the same argument as for $\Mfor$, we conclude also that $F_H$ is continuous.

For $(2)$, note that $C^*_H(i)=F_H(i)$, thus there is a condition $F\in H$ such that $F(i)=C^*_H(i)$. Hence $o^{\vec{U}}(C^*_H(i))=f(i)$ by the definition of condition in $\mathbb{M}_{f}[\vec{U}]$.

For $(3)$, as usual we note that $H$ can be defined in terms of $C^*_H$ as the filter $H_{C^*_H}$ of all the conditions $F\in \mathbb{M}_{f}[\vec{U}]$ such that for every $i\leq\lambda$,
\begin{enumerate}
    \item  If $i\in Dom(F)$, then $\kappa^{(F)}_i=C^*_H(i)$.
    \item If $i\notin Dom(F)$, then $C^*_{H}(i)\in  \underset{i\in Dom(F)}{\cup}A^{(F)}_i$.
\end{enumerate}  

$(4)$ is the standard density argument given for $\Mfor$.

As for $(5)$, note that the restriction function $\phi:\mathbb{M}_{f}[\vec{U}]\rightarrow \mathbb{M}_{f\restriction\rho}[\vec{U}]$ is a projection of forcings from  the dense subset $\{F\in \mathbb{M}_f[\vec{U}]\mid \rho\in Dom(F)\}$ onto  $\mathbb{M}_{f\restriction\rho}[\vec{U}]$, which suffices to conclude $(5)$.$\blacksquare$

The following theorem is a Mathias criteria for $\mathbb{M}_f[\vec{U}]$.
\begin{theorem}\label{MathisForMfU}
Let $f:\lambda+1\rightarrow \kappa\in V$ be a derived suitable function, and let $C\subseteq\kappa$ be such that:
\begin{enumerate}
    \item $\otp(C)=\lambda+1$ and $C$ is continuous.
    \item For every $i\leq\lambda$, $o^{\vec{U}}(C(i))=f(i)$.
    \item For every $\delta\in \Lim(\lambda+1)$, and $A\in\cap\vec{U}(C(\delta))$, there is $\xi<\delta$ such that $C\cap(\xi,\delta)\subseteq A$.
\end{enumerate}
Then there is a $V$-generic filter $H\subseteq\mathbb{M}_f[\vec{U}]$ such that $C^*_H=C$.
\end{theorem}

\pr

Define $H_C$ to consist of all the conditions $F\in \mathbb{M}_f[\vec{U}]$ such that for every $i\in Dom(F)$:
\begin{enumerate}
    \item $F(i)=C(i)$.
    \item $C\setminus \{\kappa^{(F)}_i\mid i\in Dom(F)\}\subseteq \underset{i\in Dom(F)}{\bigcup}A^{(F)}_i$.
\end{enumerate}
We prove by induction on $\lambda$ that $H_C$ is $V$-generic.  Assume for every $\rho<\lambda$ and any suitable function $g:\rho+1\rightarrow \kappa$, every $C'$ satisfying $(1)-(3)$, the definition of $H_{C'}$ is generic for $\mathbb{M}_g[\vec{U}]$. Let $f,C$ as in the theorem.
For every $\delta<\lambda$, by definition, $H_C\restriction\delta+1=H_{C\restriction\delta+1}$. Hence by the induction hypothesis $H_C\restriction\delta+1$ is generic for $\mathbb{M}_{f\restriction \delta+1}[\vec{U}]$.
 Also it is a straightforward verification that $H_C$ is a filter. 
 Let $D$ be a dense open subset of  $\mathbb{M}_f[\vec{U}]$. 
\begin{claim}
For every $F\in \mathbb{M}_f[\vec{U}]$, there is $F\leq G_F$ such that
\begin{enumerate}
    \item $\xi:=\max(Dom(F)\cap\lambda))=\max(Dom(G_F)\cap\lambda)$.
    \item There are $\xi<i_1<...<i_k<\lambda+1$ such that every $\vec{\alpha}\in\prod_{j=1}^k A^{(F)}_{\lambda,f(i_j)}$,  $G_F^{\smallfrown}\vec{\alpha}\in D$.
\end{enumerate}
\end{claim}
\pr For every $i_1<...<i_k<\lambda+1$ and every $F\leq G$ such that $$\max(Dom(F)\cap\lambda)=\max(Dom(G)\cap\lambda)\text{ and }G(\lambda)=F(\lambda)$$ consider the set
$$B=\{\vec{\alpha}\in \prod_{j=1}^k A^{(G)}_{\lambda,f(i_j)}\mid\exists R. G^{\smallfrown}\vec{\alpha}\leq^*R\in D\}$$
Then $$B\in \prod_{j=1}^k U(\kappa^{(F)}_\lambda,f(i_j))\ \ \vee\ \ \prod_{j=1}^k A^{(F)}_{\lambda,f(i_j)}\setminus B\in\prod_{j=1}^k U(\kappa^{(F)}_\lambda,f(i_j))$$ Denote the  set which is in $\prod_{j=1}^k U(\kappa^{(F)}_\lambda,f(i_j))$ by $B'$. By normality, there are $B_{i_j}\in U(\kappa^{(F)}_\lambda,f(i_j))$ such that $\prod_{j=1}^kB_{i_j}\subseteq B'$. Let $A^*_{G,i_1,..,i_k}\in\cap\vec{U}(\kappa^{(F)}_\lambda)$ be the set obtained by shrinking only the sets $A^{(F)}_{\lambda,f(i_j)}$ to $B_{i_j}$. Since $o^{\vec{U}}(\kappa^{(F)}_\lambda)<\kappa^{(F)}_\lambda$ the possibilities for $G$ (note that $G(\lambda)$ must be $F(\lambda)$) and $i_1,...,i_k$ are at most $\lambda$. So by $\kappa^{(F)}_\lambda$-completness $$A^*=\cap_{G,i_1,..,i_k}A^*_{G,i_1,...,i_k}\in \cap\vec{U}(\kappa^{(F)}_\lambda)$$ Let $F\leq^* F^*$ be the condition obtained by shrinking $A^{(F)}_\lambda$ to $A^*$. By density, there is $G\geq F$ such that $G\in D$. So there is $\vec{\alpha}\in [A^*]^{<\omega}$ such that $$\Big(G\restriction \max(Dom(F)\cap\lambda)\Big)\cup \{\l\lambda,\l \kappa^{(F)}_\lambda,A^*\r\}^{\smallfrown}\vec{\alpha}\leq^* G$$ Let $i_j\in Dom(G)$ be such that $\kappa^{(G)}_{i_j}=\vec{\alpha}(j)$, then $o^{\vec{U}}(\alpha_j)=f(i_j)$ and $\vec{\alpha}\in \prod_{j=1}^kA^{(F^*)}_{\lambda,f(i_j)}$. Hence for every $\vec{\beta}\in \prod_{j=1}^kA^{(F^*)}_{\lambda,f(i_j)}$, there is $G_{\vec{\beta}}$ such that $$\Big(G\restriction \max(Dom(F)\cap\lambda)\Big)\cup \{\l\lambda,\l \kappa^{(F)}_\lambda,A^*\r\}^{\smallfrown}\vec{\beta}\leq^*G_{\vec{\beta}}\in D$$ Note that $\vec{\beta}\in [A^*]^{<\omega}$, hence we are in the same situation as in proposition \ref{amalgamate1}, so we can find a single $F\leq G_F$ as wanted.$\blacksquare$

For every possible lower part $F_0$ below $C(\lambda)$ i.e. $F_0=F\restriction \lambda$ for some $F\in \mathbb{M}_f[\vec{U}]$ with $\kappa^{(F)}_\lambda=C(\lambda)$, use the claim to find $F_0\cup\{\l\lambda, \l C(\lambda),C(\lambda)\r\r\}\leq G_{F_0}$.
Let $$A^*=\Delta_{F_0}A_{F_0}:=\Big\{\alpha< C(\lambda)\mid \forall F_0. F_0(\max(Dom(F_0)))<\alpha\rightarrow \alpha\in A_{F_0}\Big\}\in \cap\vec{U}(C(\lambda))$$ There is $
\xi<C(\lambda)$ such that $C\cap(\xi,C(\lambda))\subseteq A^*$. Pick any $\kappa'\in C\cap[\xi,C(\lambda))$ and let $\delta<\lambda$ be such that $C(\delta)=\kappa'$. 
By the claim, the set
$$E=\Big\{F\in\mathbb{M}_{f\restriction\delta+1}[\vec{U}]\mid \exists \delta<i_1<...<i_k. \ \forall \vec{\alpha}\in\prod_{j=1}^k A^{*}_{f(i_j)}. \ G_F^{\smallfrown}\vec{\alpha}\in D\Big\}$$
is dense.
Since $H_C\restriction\delta+1$ is generic, there is $G^*\in (H_C\restriction \xi+1)\cap E$.  By condition $(2)$ of the assumption of the theorem, $f(i_j)=o^{\vec{U}}(C(i_j))$ and since $\xi<i_1<...<i_k$,  $\l C(i_1),C(i_2),...,C(i_k)\r\in \prod_{j=1}^k A^*_{f(i_j)}$. Thus $$(G^*\cup\{\langle \lambda,\l \kappa, A^*\r\r\})^{\smallfrown}\langle C(i_1),C(i_2),...,C(i_k)\rangle\in H_C\cap D$$
which conclude the proof that $H_C$ is generic. Obviously condition $(1)$ of the definition of $H_C$ ensures that $C^*_{H_C}=C$. $\blacksquare$

\begin{theorem}\label{genericover}
Let $G\subseteq \Mfor$ be $V$-generic and let $C\subseteq C_G$ be any closed subset. Let $f_C$ be the suitable function derived from $C$. If $f_C\in V$, then there is a $V$-generic $H\subseteq\mathbb{M}_{f_C}[\vec{U}]$ such that $C^*_H=C$.
\end{theorem}
\pr Let us certify that $C$ satisfy the assumptions of theorem \ref{MathisForMfU} with respect to $f_C$. $(1),(2)$ are immediate from the definition of $f_C$ and by closure of $C$. To see condition $(3)$, let $\delta\in \Lim(\lambda+1)$ and $A\in\cap\vec{U}(C(\delta))$. Since $C(\delta)\in \Lim(C)$, and $C\subseteq C_G$, $C(\delta)\in \Lim(C_G)$. By proposition \ref{genericproperties}(3), there is $\xi<\delta$ such that $C_G\cap(\xi,\delta)\subseteq A$ and also $C\cap(\xi,\delta)\subseteq A$. $\blacksquare$

\begin{example}
Consider the Prikry forcing with $U(\kappa,0)$, take $C=C_G\restriction_{even}$. Then $\otp(C\cup\{\kappa\})=\omega+1$, $f_{C}(n)=o^{\vec{U}}(C_G(2n))=0$ and $f_C(\omega)=o^{\vec{U}}(\kappa)>0$. The forcing $\mathbb{M}_{f_C}[\vec{U}]$ is simply the Prikry forcing with $U(\kappa,0)$.  Distinguishing from the forcing $\MforI$, where we must leave ``room" for the missing elements of the full generic $C_G$, it is possible that $\mathbb{M}_{f_C}[\vec{U}]$ did not leave ordinals between successive points of the Prikry sequence.
\end{example}
\begin{theorem}\label{theorem of classification}
Assume that $\forall\alpha\leq\kappa.o^{\vec{U}}(\alpha)<\alpha$. Let $G\subseteq\Mfor$ be $V$-generic filter and let $V\subseteq M\subseteq V[G]$ be an intermediate $ZFC$ model. Then there is a closed subset $C^*_{fin}\subseteq C_G$ such that $M=V[C^*_{fin}]$ and $V[C^*_{fin}]$ is a generic extension of a finite iteration of the form  
$$\mathbb{M}_{f_1}[\vec{U}]*\mathbb{M}_{\underaccent{\sim}{f}_2}[\vec{U}]...*\mathbb{M}_{\underaccent{\sim}{f}_n}[\vec{U}]$$
\end{theorem}
\pr
By \cite[Thm. 15.43]{Jech2003}, there is $A\in V[G]$ such that $V[A]=M$. By theorem \ref{MainResaultPartone}, there is $C'\subseteq C_G$ such that $M=V[A]=V[C']$.
Apply \ref{FindingClosedSet} to find a closed $C^*\subseteq C_G\cup\{\kappa\}$ such that $V[C']=V[C^*]$. Let $\lambda_0=\kappa$, recursively define
$\lambda_{i+1}=\otp(C_G\cap \lambda_i)$.
By the assumption $\forall\alpha\leq\kappa.o^{\vec{U}}(\alpha)<\alpha$ and proposition \ref{indc}, $\otp(C_G\cap\lambda_i)<\lambda_i$. Hence
after finitely many steps, $\lambda_n\leq C_G(\omega)$, denote $\kappa_i=\lambda_{n-i}$. Let $C^*_n:=C^*$ and
consider the derived suitable function $$f_n:=f_{C^*_n\cap(\kappa_{n-1},\kappa_n]}:\otp(C^*_n\cap(\kappa_{n-1},\kappa_n])\rightarrow \kappa$$
Since for each $x\in C^*_n\cap(\kappa_{n-1},\kappa_n)$, $o^{\vec{U}}(x)<\otp(C_G\cap\kappa_n)$ and $\otp(C^*\cap(\kappa_{n-1},\kappa_n))\leq \kappa_{n-1}$. By proposition \ref{genericproperties}(6), $f_n\in V[C^*_n]\cap V[C_G\cap \kappa_{n-1}]$. 
By proposition \ref{MainResaultPartone} there is $D\subseteq C_G\cap \kappa_{n-1}$ such that $V[f_n]=V[D]$, apply proposition \ref{union of generic} to $D,C^*_n\cap\kappa_{n-1}$  to find $E\subseteq \kappa_{n-1}$ such that $V[D,C^*_n\cap\kappa_{n-1}]=V[E]$. Next, apply \ref{FindingClosedSet} to $E$ in order
to find a closed subset $C^*_{n-1}\subseteq C_G\cap \kappa_{n-1}\cup\{\kappa\}$ such that $V[C^*_{n-1}]=V[E]$. Now consider the derived suitable function $$f_{n-1}:=f_{C^*_{n-1}\cap(\kappa_{n-2},\kappa_{n-1}]}:\otp(C^*_{n-1}\cap(\kappa_{n-2},\kappa_n-1])\rightarrow \kappa$$ By the same arguments as before, $f_{n-1}\in V[C^*_{n-1}]\cap V[C_G\cap\kappa_{n-2}]$ and there is a closed subset $C^*_{n-2}\subseteq C_G\cap \kappa_{n-2}\cup\{\kappa_{n-2}\}$ such that $C^*_{n-2}\in V[C^*_{n-1}]$ and $V[C^*_{n-2}]=V[C^*_{n-1}\cap \kappa_{n-2},f_{n-1}]$.
In a similar fashion we define $C^*_0,C^*_1,...,C^*_n$ such that:
\begin{enumerate}
    \item  For every $0\leq i\leq n$, $C^*_i\subseteq C_G\cap\kappa_i\cup\{\kappa_i\}$ is closed.
    \item $V[C^*_0]\subseteq V[C^*_1]\subseteq V[C^*_2]...\subseteq V[C^*_n]=M$.
    \item For every $0\leq i\leq n$, $V[C^*_{i}]=V[C^*_{i+1}\cap\kappa_i,f_{i+1}]$, where $f_{i+1}=f_{C^*_{i+1}\cap(\kappa_{i},\kappa_{i+1}]}$.
    \item 
    $f_0\in V$.
\end{enumerate}
Item $(4)$ follows from $C^*_0\subseteq\{ C_G(n)\mid n<\omega\}$.
$$C^*_{fin}=C^*_0\uplus( C^*_1\setminus \kappa_0)\uplus( C^*_2\setminus \kappa_1)\uplus....\uplus(C^*_n\setminus \kappa_{n-1})$$
\begin{claim}

\begin{enumerate}
    \item $C^*_{fin}$ is closed.
    \item For every $0\leq i\leq n$, $V[C^*_{fin}\cap\kappa_i]=V[C^*_i]$ and in particular $V[C^*_{fin}]=V[C^*]=M$
    \item For every $0<i\leq n$, $f_i=f_{C^*_{fin}\cap(\kappa_{i-1},\kappa_i]}\in V[C^*_{fin}\cap\kappa_{i-1}]$.
    
\end{enumerate}
\end{claim}
\pr $C^*_{fin}$ is a closed as the union of finitely many close sets. We prove $(2)$ by induction, for $i=0$, $C^*_{fin}\cap \kappa_0=C^*_0$.
Assume that $V[C^*_{fin}\cap \kappa_i]=V[C^*_i]$, 
then $$V[C^*_{fin}\cap\kappa_{i+1}]=V[C^*_{fin}\cap \kappa_i,C^*_{fin}\cap (\kappa_i,\kappa_{i+1})]=V[C^*_i, C^*_{i+1}\setminus \kappa_i].$$
To see that $V[C^*_i, C^*_{i+1}\setminus \kappa_i]=V[C^*_{i+1}]$, we use the third property of the sequence $C^*_j$, namely that $V[C^*_{i}]=V[C^*_{i+1}\cap\kappa_i,f_{i+1}]$ to see that 
$C^*_{i+1}\in V[C^*_i, C^*_{i+1}\setminus \kappa_i]$ and therefore $C^*_{i+1}\in V[C^*_i, C^*_{i+1}\setminus \kappa_i]$. As for the other direction, by the second property, $C^*_i\in V[C^*_{i+1}]$ and also $C^*_{i+1}\setminus \kappa_i\in V[C^*_{i+1}]$, so we conclude that $V[C^*_{fin}\cap\kappa_{i+1}]= V[C^*_{i+1}]$.

As for $(3)$, note that $C^*_{fin}\cap(\kappa_{i-1},\kappa_i]=C^*_i\cap(\kappa_{i-1},\kappa_i]$, and by property $(3)$ of the sequence $C^*_j$, $f_i\in V[C^*_{i-1}]$. By $(2)$ of the claim it follow that $$f_{C^*_{fin}\cap(\kappa_{i-1},\kappa_i]}=f_{C^*_i\cap(\kappa_{i-1},\kappa_i]}=f_i\in V[C^*_{i-1}]=V[C^*_{fin}\cap\kappa_{i-1}]\ \ \ \ \ \ \ \ \ \ \ \ \ \blacksquare_{claim}$$
 Therefore for every $i\leq n$, $\mathbb{M}_{f_i}[\vec{U}]$ is defined in $V[C^*_{fin}\cap\kappa_{i-1}]$, denote this model by $N_i$. Recall remark \ref{lift}, the club $C_G\cap (\kappa_{i-1},\kappa_i)$ is $V[C_G\cap\kappa_{i-1}]$-generic for the forcing $\Mfor\restriction(\kappa_{i-1},\kappa_i)$\footnote{ Alternatively, it is $V[C_G\cap\kappa_{i-1}]$-generic for  $\mathbb{M}[\vec{W}]\restriction(\kappa_{i-1},\kappa_i)$, where $\vec{W}$ is the coherent sequence generated by $\vec{U}$ in $V[C_G\cap\kappa_{i-1}]$.}  and therefore it is $N_i$-generic as $N_i\subseteq V[C_G\cap\kappa_{i-1}]$. Hence we can
 apply theorem \ref{genericover} to $C^*_{fin}\cap(\kappa_{i-1},\kappa_i]\subseteq C_G\cap(\kappa_{i}+1)$ and find a
 $N_i$-generic filter $H\subseteq\mathbb{M}_{f_i}[\vec{U}]$ such that $$N_i[H]=N_i[C^*_{fin}\cap(\kappa_{i-1},\kappa_i]=V[C^*_{fin}\cap\kappa_{i-1}][C^*_{fin}\cap(\kappa_{i-1},\kappa_i]]=V[C^*_{fin}\cap\kappa_i]$$
 In particular, $V[C^*_{fin}\cap\kappa_0]$ is a generic extension of $V$ by $\mathbb{M}_{f_0}[\vec{U}]$.

 Let $\underaccent{\sim}{f}_i$ be a $(\mathbb{M}_{f_0}[\vec{U}]*\mathbb{M}_{\underaccent{\sim}{f}_1}[\vec{U}]...*\mathbb{M}_{\underaccent{\sim}{f}_{i-1}}[\vec{U}])$-name for $f_i$, then there is a $V$-generic filter $H^*$ for the iteration $\mathbb{M}_{f_1}[\vec{U}]*\mathbb{M}_{\underaccent{\sim}{f}_2}[\vec{U}]...*\mathbb{M}_{\underaccent{\sim}{f}_n}[\vec{U}]$ such that $V[H^*]=V[C^*_{fin}]=M$ (see for example \cite[Thm. 16.2]{Jech2003}).$\blacksquare_{theorem \ \ref{theorem of classification}}$
\section{Acknowledgment}

The authors would like to thank the referee for his careful examination of the paper and for many useful insights. Also they would like to thank the participants of the ``Set Theory 
Seminar" of Tel-Aviv University for their comments during the presentation of this work. 

\bibliographystyle{amsplain}
\bibliography{ref}

\providecommand{\bysame}{\leavevmode\hbox to3em{\hrulefill}\thinspace}
\providecommand{\MR}{\relax\ifhmode\unskip\space\fi MR }
\providecommand{\MRhref}[2]{%
  \href{http://www.ams.org/mathscinet-getitem?mr=#1}{#2}
}
\providecommand{\href}[2]{#2}
\begin{thebibliography}{1}

\bibitem{TomMoti}
Tom Benhamou and Moti Gitik, \emph{{S}ets in {P}rikry and {M}agidor {G}eneric
  {E}xtesions}, Annals of Pure and Applied Logic \textbf{172} (2021), no.~4,
  102926.

\bibitem{Gitik2010}
Moti Gitik, \emph{{P}rikry-{T}ype {F}orcings}, pp.~1351--1447, Springer
  Netherlands, Dordrecht, 2010.

\bibitem{PrikryCaseGitikKanKoe}
Moti Gitik, Vladimir Kanovei, and Peter Koepke, \emph{{I}ntermediate {M}odels
  of {P}rikry {G}eneric {E}xtensions}, Pre Print (2010),
  http://www.math.tau.ac.il/~gitik/spr--kn.pdf.

\bibitem{Jech2003}
Thomas Jech, \emph{Set {T}heory}, Springer Monographs in Mathematics,
  Springer-Verlag, Berlin, 2003, The third millennium edition, revised and
  expanded. \MR{1940513}

\bibitem{ChangeCofinality}
Menachem Magidor, \emph{{C}hanging the {C}ofinality of {C}ardinals}, Fundamenta
  Mathematicae \textbf{99} (1978), 61--71.

\bibitem{MitchellMathias}
William Mitchell, \emph{{H}ow {W}eak is a {C}losed {U}nbounded {F}ilter?},
  stud. logic foundation math. \textbf{108} (1982), 209--230.

\bibitem{Unger2013}
Spencer Unger, \emph{Aronszajn trees and the successors of a singular
  cardinal}, Arch. Math. Logic \textbf{52} (2013), no.~5-6, 483--496.
  \MR{3072773}

\end{thebibliography}

\end{document}